\documentclass[11pt]{article}
\usepackage{amsmath,amsthm,amsfonts,amssymb,latexsym}
\usepackage{amsmath}

\usepackage{a4wide}
\newtheorem{thm}{Theorem}[section]

\newtheorem{lem}[thm]{Lemma}
\newtheorem{prop}[thm]{Proposition}

\newtheorem{theorem}{Theorem}[section]
\newtheorem{definition}[theorem]{Definition}
\newtheorem{lemma}[theorem]{Lemma}
\newtheorem{remark}[theorem]{Remark}
\newtheorem{proposition}[theorem]{Proposition}
\newtheorem{corollary}[theorem]{Corollary}

\def\a{\mathfrak a}
\def\k{\mathfrak k}

\def\a{\alpha}

\def\l{\lambda}

\def\k{\mathfrak{k}}

\def\L{{\cal   L}}

\def\HH{{\rm   H}}
\def\fF{{\frak   F}}

\def\Vect{{\rm   Vect}}
\def\cK{{\cal   K}}

\def\half{{\frac{1}{2} }}
\def\cF{{\cal   F}}

\def\L1#1{L^1(#1)}

\def\lef({\left(}
\def\rig){\right)}

\begin{document}

\title{  Cohomology of $\frak {osp}(1|2)$ acting
on the space of bilinear differential operators on the superspace
$\mathbb{R}^{1|1}$ }

\label{firstpage}

\author{ Mabrouk Ben Ammar \and Amina Jabeur \and Imen Safi \thanks{D\'epartement de Math\'ematiques,
Facult\'e des Sciences de Sfax, BP 802, 3038 Sfax, Tunisie.
E.mail:  mabrouk.benammar@fss.rnu.tn}}

%--------------------------------------------------------------------
\maketitle
% ----------------------------------------------------------------

\begin{abstract}
We compute the first cohomology of the
ortosymplectic Lie superalgebra $\mathfrak{osp}(1|2)$ on the
(1,1)-dimensional real superspace with coefficients in the superspace
$\frak{D}_{\lambda,\nu;\mu}$ of bilinear differential operators
acting on weighted densities. This work is the simplest
superization of a result by Bouarroudj [Cohomology of the vector fields Lie algebras on $\mathbb{R}\mathbb{P}^1$ acting on bilinear differential operators, International Journal of Geometric Methods in Modern Physics
 (2005), {\bf 2}; N 1,  23-40].

\end{abstract}

\maketitle {\bf Mathematics Subject Classification} (2000). 53D55


\thispagestyle{empty}

%%%%%%%%%%%%%%%%%%%%%%%%%%%%%%%%%%%%%%%%%%%%%%%%%%%%%%%%%%%%%%%%%%%%%%%%%%%%%%%
%%%%%%%%%%%%%%%%%%%%%%%%%%%%%%%%%%%%%%%%%%%%%%%%%%%%%%%%%%%%%%%%%%%%%%%%%%%%%%%
\section{Introduction}
%%%%%%%%%%%%%%%%%%%%%%%%%%%%%%%%%%%%%%%%%%%%%%%%%%%%%%%%%%%%%%%%%%%%%%%%%%%%%%%
%%%%%%%%%%%%%%%%%%%%%%%%%%%%%%%%%%%%%%%%%%%%%%%%%%%%%%%%%%%%%%%%%%%%%%%%%%%%%%%

The space of weighted densities of weight $\lambda$ on
$\mathbb{R}$ (or $\lambda$-densities for short), denoted by:
\begin{equation*}{\cal F}_\lambda=\left\{ fdx^{\lambda}, ~f\in
C^\infty(\mathbb{R})\right\},\quad \lambda\in\mathbb{R},
\end{equation*}is the space of
sections of the line bundle $(T^*\mathbb{R})^{\otimes^\l}.$ The
Lie algebra ${\rm Vect}(\mathbb{R})$ of vector fields $X_h=h{d\over
dx}$, where $h\in C^\infty(\mathbb{R})$, acts by the {\it Lie
derivative}. Alternatively, this action can be written as follows:
\begin{equation}\label{Lie1}{X_h}\cdot(fdx^{\lambda})=L_{X_h}^\lambda(f)dx^{\lambda}\quad\text{with}~~L_{X_h}^\lambda(f)=hf'+\lambda
h'f,
\end{equation}
where $f'$, $h'$ are $\frac{df}{dx}$, $\frac{dh}{dx}$. Each bilinear
differential operator $A$ on $\mathbb{R}$ gives thus rise to a
morphism from $\cF_\l\otimes\cF_\nu$ to $\cF_\mu$, for any
$\l,\,\nu,\,\mu\in\mathbb{R}$, by $fdx^\lambda\otimes gdx^\nu\mapsto
A(f\otimes g)dx^\mu$. The Lie algebra ${\rm Vect}(\mathbb{R})$ acts
on the space $\mathrm{D}_{\lambda,\nu;\mu}$ of these differential
operators by:
\begin{equation}\label{Lieder2}X_h\cdot A=L_{X_h}^\mu\circ
A-A\circ L_{X_h}^{(\lambda,\nu)} \end{equation} where
$L_{X_h}^{(\lambda,\nu)}$ is the Lie derivative on
$\cF_\l\otimes\cF_\nu$  defined by the Leibnitz rule:
$$L_{X_h}^{(\lambda,\nu)}(f\otimes
g)=L_{X_h}^{\lambda}(f)\otimes
g+f\otimes L_{X_h}^{\nu}(g).$$
If we restrict ourselves to the Lie algebra $\frak{ sl}(2)$ which
is isomorphic to the Lie subalgebra of ${\rm Vect}(\mathbb{R})$
spanned by
\begin{equation*}\left\{X_1,\,X_x,\,X_{x^2}\right\},\end{equation*}
we have a family of infinite dimensional $\frak{ sl}(2)$-modules
still denoted by $\mathrm{D}_{\lambda,\nu;\mu}$. Bouarroudj, in
\cite {b}, computes the cohomology space
$\mathrm{H}^1_\mathrm{diff}\left(\mathfrak{sl}(2),
\mathrm{D}_{\lambda,\nu;\mu}\right)$ where
$\mathrm{H}^1_\mathrm{diff}$ denotes the differential cohomology;
that is, only cochains given by differential operators are
considered.

In this paper we are interested in the study of the analogue super
structures. More precisely, we consider here the superspace
$\mathbb{R}^{1|1}$ equipped with the standard {\it contact
structure} given by the 1-form $\alpha=dx+\theta d\theta$, we
replace $\frak{ sl}(2)$ by its analogue in the super setting, i.e
the orthosymplectic  Lie superalgebra $\mathfrak{osp}(1|2)$ which
can be realized as a subalgebra of the superalgebra $\mathcal{K}(1)$
of contact vector fields. We introduce the superspace of
$\lambda$-densities on the superspace $\mathbb{R}^{1|1}$ (with
respect to $\alpha$) denoted by $\fF_{\l}$ and the superspace
$\frak{D}_{\lambda,\nu;\mu} $ of differential bilinear operators
viewed as homomorphisms from $\fF_{\l}\otimes\fF_\nu$ to
$\fF_{\mu}$. The superalgebra $\mathfrak{osp}(1|2)$   acts naturally
on $\fF_{\l}$ and $\frak{D}_{\lambda,\nu;\mu} $. We compute here the
first cohomology spaces $\mathrm{H}^1_\mathrm{diff}\left(\frak
{osp}(1|2) ,\mathfrak{sl}(2);\frak{D}_{\l,\nu;\mu}\right) $ and
$\mathrm{H}^1_\mathrm{diff}\left(\frak {osp}(1|2)
,\frak{D}_{\l,\nu;\mu}\right) $, $\l,\,\nu,\,\mu\in\mathbb{R}$,
getting a result very close to the classical one
$\mathrm{H}^1_\mathrm{diff}\left(\mathfrak{sl}(2),
\mathrm{D}_{\lambda,\nu;\mu}\right)$. Moreover, we give explicit
formulae for non trivial 1-cocycles which generate these spaces.

These spaces appear naturally in the problem of describing the
deformations of the $\mathfrak{osp}(1|2)$-modules
$\frak{D}_{\l,\nu;\mu}$. More precisely, the first cohomology space
$\mathrm{H}^1\left(\mathfrak{osp}(1|2),V\right)$ classifies the
infinitesimal deformations of an $\mathfrak{osp}(1|2)$ module $V$
and the obstructions to integrability of a given infinitesimal
deformation of $V$ are elements of
$\mathrm{H}^2\left(\mathfrak{osp}(1|2),V\right)$.

 %%%%%%%%%%%%%%%%%%%%%%%%%%%%%%%%%%%%%%%%%%%%%%%%%%%%%%%%%%%%%%%%%%
\section{Definitions and Notation}
%%%%%%%%%%%%%%%%%%%%%%%%%%%%%%%%%%%%%%%%%%%%%%%%%%%%%%%%%%
%%%%%%%%%%%%%%%%%%%%%%%%%%%%%%%%%%%%%%%%%%%%%%%%%%%%%%%%%%
\subsection{The Lie superalgebra of contact vector fields on $\mathbb{R}^{1|1}$}
%%%%%%%%%%%%%%%%%%%%%%%%%%%%%%%%%%%%%%%%%%%%%%%%%%%%%%%%%%
We define the superspace $\mathbb{R}^{1|1}$ in terms of its
superalgebra of functions, denoted by $C^\infty(\mathbb{R}^{1|1})$
and consisting of elements of the form:
\begin{equation*}
F(x,\theta)=f_0(x)+f_1(x)\theta ,
\end{equation*}
where $x$ is the even variable,  $\theta$ is  the odd variable
($\theta^2=0$) and $f_0(x),\,f_1(x)\in C^\infty(\mathbb{R})$. Even elements in $C^\infty(\mathbb{R}^{1|1})$ are the functions
$F(x,\theta)=f_0(x)$, the functions $F(x,\theta)=\theta f_1(x)$ are odd elements. The parity of homogenous elements $F$ will be denoted $|F|$. We
consider the contact bracket on $C^\infty(\mathbb{R}^{1|1})$
defined on  $C^\infty(\mathbb{R}^{1|1})$ by:
\begin{equation*}\begin{array}{l}
\{F,G\}=FG'-F'G+\frac{1}{2}\eta(F)
\overline{\eta}(G),\end{array}
\end{equation*}where $\eta=\frac{\partial}{\partial
{\theta}}+\theta\frac{\partial}{\partial x}$ and
$\overline{\eta}=\frac{\partial}{\partial
{\theta}}-\theta\frac{\partial}{\partial x}$. The superspace $\mathbb{R}^{1|1}$ is equipped with
the standard contact structure given by the following $1$-form:
\begin{equation*}
\a=dx+\theta d\theta.
\end{equation*}
Let $\mathrm{Vect}(\mathbb{R}^{1|1})$ be the superspace of
 vector fields on $\mathbb{R}^{1|1}$:
\begin{equation*}\mathrm{Vect}(\mathbb{R}^{1|1})=\left\{F_0\partial_x
+  F_1\partial_\theta \mid ~F_i\in C^\infty(
\mathbb{R}^{1|1})\right\},\end{equation*} where $\partial_\theta$
stands for $\frac{\partial}{\partial\theta}$ and $\partial_x$
stands for $\frac{\partial}{\partial x} $, and consider the
superspace $\mathcal{K}(1)$ of contact vector fields on
$\mathbb{R}^{1|1}$. That is, $\mathcal{K}(1)$ is the Lie
superalgebra of conformal vector fields on $\mathbb{R}^{1|1}$ with
respect to the $1$-form $\alpha$:
 $$\mathcal{K}(1)=\big\{X\in\mathrm{Vect}(\mathbb{R}^{1|1})~|~\hbox{there
exists}~H\in C^\infty(\mathbb{R}^{1|1})~ \hbox{such
that}~\mathfrak{L}_X(\alpha)=H\alpha\big\},$$
where $\mathfrak{L}_X$ is the Lie derivative along the vector field $X$.
Any contact vector field
on $\mathbb{R}^{1|1}$ has the following explicit form:
\begin{equation*}\begin{array}{l}
X_H=H\partial_x+\half\eta(H)\overline{\eta},\;\text{ where }\, H\in
C^\infty(\mathbb{R}^{1|1}).\end{array}
\end{equation*}
The bracket on $\mathcal{K}(1)$ is given by
\begin{equation*}
[X_{F},\,X_{G}]=X_{\{F,\,G\}}.
\end{equation*}

%%%%%%%%%%%%%%%%%%%%%%%%%%%%%%%%%%%%%%%%%%%%%%%%%%%%%%%%%%%%%%%%%%%%%%%%%%%%%%%%%%%%%%%%%%%%%
\subsection{The subalgebra $\mathfrak{osp}(1|2)$}
%%%%%%%%%%%%%%%%%%%%%%%%%%%%%%%%%%%%%%%%%%%%%%%%%%%%%%%%%%%%%%%%%%%%%%%%%%%%%%%%%%%%%%%%%%%%%

The Lie algebra $\frak {sl}(2)$ is realized as subalgebra of the
Lie algebra  $\Vect(\mathbb{R})$:
\begin{equation*}\frak {sl}(2)=\text{Span}(X_1,\,X_{x},\,X_{x^2}).\end{equation*}
Similarly, we now consider the orthosymplectic Lie superalgebra as
a subalgebra of $\mathcal{K}(1)$:
\begin{equation*}
\mathfrak{osp}(1|2)=\text{Span}(X_1,\,X_{x},\,X_{x^2},\,X_{x\theta},\,
X_{\theta}).
\end{equation*}
The space of even elements is isomorphic to $\mathfrak{sl}(2)$, while the space of odd elements is two dimensional:
\begin{equation*}
(\mathfrak{osp}(1|2))_{\bar1}=\text{Span}(X_{x\theta},\,X_{\theta}).
\end{equation*}
The new commutation relations are
\begin{equation*}
\begin{array}{llll}
&[X_{x^2},X_\theta]~=-X_{x\theta},~~&[X_x,X_\theta]~=-\frac{1}{2}
X_\theta,~~&[X_1,X_\theta]=0,\cr
&[X_{x^2},X_{x\theta}]=0,~~&[X_x,X_{x\theta}]=\frac{1}{2}
X_{x\theta},~~&[X_1,X_{x\theta}]=X_\theta,\cr
&[X_{x\theta},X_{\theta}]~=\frac{1}{2}
X_{x},~~&[X_{\theta},X_{\theta}]~=\frac{1}{2}
X_{1},~~&[X_{x\theta},X_{x\theta}]=\frac{1}{2} X_{x^2}.
\end{array}
\end{equation*}

%%%%%%%%%%%%%%%%%%%%%%%%%%%%%%%%%%%%%%%%%%%%%%%%%%%%%%%%%%%%%%%%%%%%%%%%%%%%%
\subsection{The space of weighted densities on $\mathbb{R}^{1|1}$}
%%%%%%%%%%%%%%%%%%%%%%%%%%%%%%%%%%%%%%%%%%%%%%%%%%%%%%%%%%%%%%%%%%%%%%%%%%%%%
%%%%%%%%%%%%%%%%%%%%%%%%%%%%%%%%%%%%%%%%%%%%%%%%%%%%%%%%%%%%%%%%%%%%%%%%%%%%%%%

We have analogous definition of weighted
densities in super setting (see \cite{ab}) with $dx$ replaced by
$\alpha$. The elements of these spaces are indeed (weighted)
densities since all spaces of generalized tensor fields have just
one parameter relative $\cK(1)$
--- the value of $X_x$ on the lowest weight vector (the one
annihilated by $X_\theta$). From this point of view the volume
element (roughly speaking, $\lq\lq
dx\frac{\partial}{\partial\theta}"$) is indistinguishable from
$\a^{\frac{1}{2}}.$ We
denote by $\fF_{\l}$ the space of all
 weighted densities on $\mathbb{R}^{1|1}$ of weight $\lambda$:
\begin{equation*}
\mathfrak{F}_\lambda=\left\{F(x,\theta)\alpha^\lambda~~|~~F(x,\theta)
\in C^\infty(\mathbb{R}^{1|1})\right\}.
\end{equation*}
As a vector space, $\mathfrak{F}_\lambda$ is isomorphic to
$C^\infty(\mathbb{R}^{1|1})$, but the Lie derivative of the
density $F\alpha^\lambda$ along the vector field $X_H$ in
$\mathcal{K}(1)$ is now:
\begin{equation*}
\label{superaction}
\mathfrak{L}_{X_H}(F\alpha^\lambda)=\mathfrak{L}^{\lambda}_{X_H}(F)\alpha^\lambda,
\quad\text{with}~~\mathfrak{L}^{\lambda}_{X_H}(F)=\mathfrak{L}_{X_H}(F)+
\lambda H'F.
\end{equation*}
Or, if we put $H(x,\theta)=a(x)+b(x)\theta$,
$F(x,\theta)=f_0(x)+f_1(x)\theta$,
\begin{equation*}
\label{deriv}
\mathfrak{L}^{\lambda}_{X_H}(F)=L^{\lambda}_{a\partial_x}(f_0)+\frac{1}{2}~bf_1+
\left(L^{\lambda+\frac{1}{2}}_{a\partial_x}(f_1)+\lambda
f_0b'+\frac{1}{2} f'_0 b\right)\theta.
\end{equation*}
Especially, we have
\begin{equation}\label{actiondecomp}
\left\{\begin{array}{lll}
\mathfrak{L}^\lambda_{X_a}(f_0)=L^\lambda_{a\partial_x}(f_0),&&
\mathfrak{L}^\lambda_{X_a}(f_1\theta)=\theta
L^{\lambda+\frac{1}{2}}_{a\partial_x} (f_1),\cr &\hbox{
~~and~~}&\cr \mathfrak{L}^{\lambda}_{X_{b\theta}}(f_0)=(\lambda
f_0b'+\frac{1}{2} f'_0b)
\theta&&\mathfrak{L}^{\lambda}_{X_{b\theta}}(f_1\theta)=\frac{1}{2}
bf_1.
\end{array}\right.\end{equation}
Of course, for all $\lambda$, $\mathfrak{F}_\lambda$ is a
$\mathcal{K}(1)$-module:
\begin{equation*}
[\mathfrak{L}^{\lambda}_{X_F},\mathfrak{L}^{\lambda}_{X_G}] =
\mathfrak{L}^{\lambda}_{[{X_F},\,X_G]}.
\end{equation*}
We thus obtain a one-parameter family of $\mathcal{K}(1)$-modules
on $C^\infty(\mathbb{R}^{1|1})$ still denoted by
$\mathfrak{F}_\lambda$.

%------------------------------------------------------------------------------------------
\subsection{Differential operators on weighted densities}

%----------------------------------------------------------------------------------
A differential operator on $\mathbb{R}^{1|1}$ is an operator on
$C^\infty(\mathbb{R}^{1|1})$ of the following form:
$$
A= \sum_{i=0}^\ell{a}_i (x, \theta)\partial_x^i +
\sum_{i=0}^\ell{b}_i (x,
\theta)\partial_x^i\partial_\theta.
$$
In \cite{gmo}, it is proved that any local operator $A$ on
$\mathbb{R}^{1|1}$ is in fact a differential operator.

Of course, any differential operator defines a linear mapping from
$\mathfrak{F}_\lambda$ to $\mathfrak{F}_\mu$ for any $\lambda$,
$\mu\in\mathbb{R}$, thus, the space of differential operators
becomes a family of $\mathcal{K}(1)$ and $\mathfrak{osp}(1|2)$
modules denoted $\mathfrak{D}_{\lambda,\mu}$, for the natural
action:
\begin{equation*}\label{d-action}
{X_H}\cdot A=\mathfrak{L}^\mu_{X_H}\circ A-
(-1)^{|A||H|}A\circ \mathfrak{L}^\lambda_{X_H}.
\end{equation*}
Similarly, we consider  a family of ${\rm \mathcal{K}}(1)$-modules
on the space $\frak{D}_{\lambda,\nu;\mu} $  of bilinear
differential operators: $~A: {\frak
F}_{\lambda}\otimes\frak{F}_{\nu}\longrightarrow{\frak F}_\mu$
with the $\mathcal{K}(1)$-action
\begin{equation*}X_H\cdot A={\frak L}_{X_H}^\mu\circ
A-(-1)^{|A||H|}A\circ {\frak
L}_{X_H}^{(\lambda,\nu)},\end{equation*}where ${\frak
L}_{X_H}^{(\lambda,\nu)}$ is the Lie derivative on ${\frak
F}_{\lambda}\otimes\frak{F}_{\nu}$ defined by the Leibnitz rule:
$$\begin{array}{ll}
{\frak L}_{X_H}^{(\lambda,\nu)}(F\otimes G)&= {\frak
L}_{X_H}^{\lambda}(F)\otimes G+(-1)^{|H||F|}F\otimes
{\frak L}_{X_H}^{\nu}(G).\end{array}
$$

%%%%%%%%%%%%%%%%%%%%%%%%%%%%%%%%%%%%%%%%%%%%%%%%%%%%%%%%%%%%%%%%%%%%%
%%%%%%%%%%%%%%%%%%%%%%%%%%%%%%%%%%%%%%%%%%%%%%%%%%%%%%%%%%%%%%%%%%%%%%%
\section{The space $\HH^1_\mathrm{diff}(\frak {osp}(1|2),\frak{D}_{\lambda,\nu;\mu})$}
%%%%%%%%%%%%%%%%%%%%%%%%%%%%%%%%%%%%%%%%%%%%%%%%%%%%%%%%%%%%%%
%%%%%%%%%%%%%%%%%%%%%%%%%%%%%%%%%%%%%%%%%%%%%%%%%%%%%%%%%%%%%%

%%%%%%%%%%%%%%%%%%%%%%%%%%%%%%%%%%%%%%%%%%%%%%%%%%%%%%%%%%%%%%%%%%%%%
\subsection{Cohomology }
%%%%%%%%%%%%%%%%%%%%%%%%%%%%%%%%%%%%%%%%%%%%%%%%%%%%%%%%%%%%%%%%%%%%%
We will compute the first cohomology space of $\frak {osp}(1|2)$
with coefficients in $\frak{D}_{\lambda,\nu;\mu}$. Let us first
recall some fundamental concepts from cohomology theory~(see, e.g.,
\cite{Fu}). Let $\frak{g}=\frak{g}_{\bar 0}\oplus \frak{g}_{\bar 1}$
be a Lie superalgebra acting on a superspace $V=V_{\bar 0}\oplus
V_{\bar 1}$ and let $\mathfrak{h}$ be a subalgebra of
$\mathfrak{g}$. (If $\frak{h}$ is omitted it assumed to be $\{0\}$.)
The space of $\frak h$-relative $n$-cochains of $\frak{g}$ with
values in $V$ is the $\frak{g}$-module
\begin{equation*}
C^n(\frak{g},\frak{h}; V ) := \mathrm{Hom}_{\frak
h}(\Lambda^n(\frak{g}/\frak{h});V).
\end{equation*}
The {\it coboundary operator} $ \delta_n: C^n(\frak{g},\frak{h}; V
)\longrightarrow C^{n+1}(\frak{g},\frak{h}; V )$ is a
$\frak{g}$-map satisfying $\delta_n\circ\delta_{n-1}=0$. The
kernel of $\delta_n$, denoted $Z^n(\mathfrak{g},\frak{h};V)$, is
the space of $\frak h$-relative $n$-{\it cocycles}, among them,
the elements in the range of $\delta_{n-1}$ are called $\frak
h$-relative $n$-{\it coboundaries}. We denote
$B^n(\mathfrak{g},\frak{h};V)$ the space of $n$-coboundaries.

By definition, the $n^{th}$ $\frak h$-relative  cohomolgy space is
the quotient space
\begin{equation*}
\mathrm{H}^n
(\mathfrak{g},\frak{h};V)=Z^n(\mathfrak{g},\frak{h};V)/B^n(\mathfrak{g},\frak{h};V).
\end{equation*}
We will only need the formula of $\delta_n$ (which will be simply
denoted $\delta$) in degrees 0 and 1: for $v \in
C^0(\frak{g},\,\frak{h}; V) =V^{\frak h}$,~ $\delta v(g) : =
(-1)^{p(g)p(v)}g\cdot v$, where
\begin{equation*}
V^{\frak h}=\{v\in V~\mid~h\cdot v=0\quad\text{ for all }
h\in\frak h\},
\end{equation*}
and  for  $ \Upsilon\in C^1(\frak{g}, \frak{h};V )$,
\begin{equation*}\delta(\Upsilon)(g,\,h):=
(-1)^{p(g)p(\Upsilon)}g\cdot
\Upsilon(h)-(-1)^{p(h)(p(g)+p(\Upsilon))}h\cdot
\Upsilon(g)-\Upsilon([g,~h])\quad\text{for any}\quad g,h\in
\frak{g}.
\end{equation*}

%%%%%%%%%%%%%%%%%%%%%%%%%%%%%%%%%%%%%%%%%%%%%%%%%%%%%%%%%%%%%%%%%%%%%
\subsection{The main theorem}
%%%%%%%%%%%%%%%%%%%%%%%%%%%%%%%%%%%%%%%%%%%%%%%%%%%%%%%%%%%%%%%%%%%%%
We will prove that non-zero spaces $\HH^1_\mathrm{diff}(\frak
{osp}(1|2),\frak{D}_{\l,\nu;\mu})$ only appear if
$2(\mu-\lambda-\nu)\in\mathbb{N}$. Moreover, if $\mu-\lambda-\nu$ is
integer then $\HH^1_\mathrm{diff}(\frak
{osp}(1|2),\frak{D}_{\l,\nu;\mu})$ is purely even and if
$\mu-\lambda-\nu$ is semi-integer then $\HH^1_\mathrm{diff}(\frak
{osp}(1|2),\frak{D}_{\l,\nu;\mu})$ is purely odd.
\begin{definition} ~

1)   We say that $(\l,\nu,\mu)$ is resonant if $\mu-\l-\nu-1= k$
with $k\in\mathbb{N}$, and
\begin{equation}
\begin{array}{l}(\lambda,\nu)=(-\frac{s}{2},-\frac{t}{2}),
\quad\text{ where }\quad s,\,t\in\{0,\,\dots,\,k\}\quad\text{ and }
\quad s+t\geq k.\end{array}
\end{equation}
We say that $(\l,\nu,\mu)$ is weakly resonant  if
$\mu-\l-\nu\in\mathbb{N}$ but $(\l,\nu,\mu)$ is not resonant.\\

2) We say that $(\l,\nu,\mu)$ is  super resonant if $\mu-\l-\nu-1=
k$ with $k\in\frac{1}{2}\mathbb{N}$, and
\begin{equation}
\begin{array}{l}(\lambda,\nu)=(-\frac{s}{2},-\frac{t}{2}),
\quad\text{ where }\quad s,\,t\in\{1,\,\dots,\,[k]\}\quad\text{ and
} \quad s+t\geq [k+\frac{1}{2}]+1.\end{array}
\end{equation}
We say that $(\l,\nu,\mu)$ is weakly super resonant if
$\mu-\l-\nu=k+1\in\frac{1}{2}\mathbb{N}$, and
\begin{equation}
\begin{array}{l}(\lambda,\nu)=(-\frac{s}{2},-\frac{t}{2}),
\, s,\,t\in\{0,\,\dots,\,[k]+1\}\,\Rightarrow \, s+t<[k+\frac{1}{2}].\end{array}
\end{equation}
\end{definition}
\begin{remark} The super resonance (respectively, weakly super resonance) of
$(\lambda,\nu,\mu)$ express the resonance (respectively, weakly resonance) of :
\begin{itemize}
  \item  $(\lambda,\nu,\mu)$, $(\lambda+\frac{1}{2},\nu+\frac{1}{2},\mu)$,
  $(\lambda+\frac{1}{2},\nu,\mu+\frac{1}{2})$ and $(\lambda,\nu+\frac{1}{2},\mu+\frac{1}{2})$
  if $\mu-\lambda-\nu$ is integer.
  \item  $(\lambda,\nu,\mu+\frac{1}{2})$, $(\lambda+\frac{1}{2},\nu+\frac{1}{2},\mu+\frac{1}{2})$,
  $(\lambda+\frac{1}{2},\nu,\mu)$ and $(\lambda,\nu+\frac{1}{2},\mu)$
  if $\mu-\lambda-\nu$ est semi-integer.
\end{itemize}\end{remark}

The main result in this paper is the following:
\begin{thm}
\label{th1}
\begin{equation*}
{\rm H}^1_\mathrm{diff}(\frak
{osp}(1|2),\mathfrak{D}_{\lambda,\nu;\mu})\simeq \left\{
\begin{array}{llllll}
\mathbb{R}^6 & \hbox{ if }\quad
 (\lambda,\nu,\mu)\,\text{ is super resonant, }\\[2pt]
\mathbb{R} & \hbox{ if }\quad
 (\lambda,\nu,\mu)\,\text{ is weakly super resonant.}
\end{array}
\right.
\end{equation*}

\end{thm}

The proof of Theorem \ref{th1} will be the subject of Section
\ref{s4}. Moreover, explicit
formulae for non trivial 1-cocycles  generating  the corresponding cohomology spaces will be given.
We will show that the spaces $\HH^1_{\rm diff}(\frak
{osp}(1|2),\,\frak{D}_{\lambda,\nu;\mu})$ and $\HH^1_{\rm
diff}(\frak {sl}(2),\,\mathrm{D}_{\lambda,\nu;\mu})$ are closely
related. Therefore, for comparison and to build upon, we need to
recall the description of $\HH^1_{\rm diff}(\frak
{sl}(2),\,\mathrm{D}_{\lambda,\nu;\mu})$. Of course, there are some cases
of $(\lambda,\nu,\mu)$ which are neither
super resonant nor weakly super resonant, these cases will be studied in Section \ref{s5}.

%%%%%%%%%%%%%%%%%%%%%%%%%%%%%%%%%%%%%%%%%%%%%%%%%%%%%%%%%%%%%%%%%%%%%
 \subsection{The space  $\HH^1_\mathrm{diff}(\frak
{sl}(2),\mathrm{D}_{\lambda,\nu;\mu})$}
 %%%%%%%%%% %%%%%%%%%%%%%%%%%%%%%%%%%%%%%%%%%%%%%%%%%%%%%%%%%%%%%%%%%%%
For the sake of simplicity, the elements $f dx^\lambda$ of
$\mathcal{F}_\lambda$  will be denoted  $f$. Any 1-cochain $c\in
Z^1_\mathrm{diff}(\frak {sl}(2),\mathrm{D}_{\lambda,\nu;\mu})$
should retains the following general form:
$$
c(X_h,f,\,g) = \sum_{i,j}\alpha_{i,j} h
f^{(i)} g^{(j)} + \sum_{i,j}\beta_{i,j} h'
f^{(i)} g^{(j)}+\sum_{i,j}\gamma_{i,j} h''
f^{(i)} g^{(j)}.
$$
So, for any integer $k\geq0$, we define the $k$-homogeneous component of $c$ by
$$
c_k(X_h,f,g) = \sum_{i+j=k}\alpha_{i,j} h
f^{(i)} g^{(j)} + \sum_{i+j=k-1}\beta_{i,j} h'
f^{(i)} g^{(j)}+\sum_{i+j=k-2}\gamma_{i,j} h''
f^{(i)} g^{(j)}.
$$
Of course, we suppose that $\gamma_{i,j}=0$ if $k\in\{0,\,1\}$ and $\beta_{i,j}=0$ if $k=0$. The coboundary map $\delta$ is homogeneous, therefore, we  easily  deduce the following lemma:
\begin{lem}\label{l32}
Any 1-cochain $c\in C^1_\mathrm{diff}(\frak
{sl}(2),\mathrm{D}_{\lambda,\nu;\mu})$ is a 1-cocycle if and only if
each of its homogeneous components  is a 1-cocycle.
\end{lem}\label{l33}
The following lemma gives the general form of any homogeneous 1-cocycle.
\begin{lem}\label{lem1}Up to a coboundary, any $(k+2)$-homogeneous 1-cocycle
$c\in\mathrm{Z}^1_\mathrm{diff}(\frak{sl}(2),\mathrm{D}_{\lambda,\nu;\mu})$
can be expressed as follows. For all $f\in\mathcal{F_\lambda}$,
$g\in\mathcal{F_\nu}$ and for all $X_h\in\frak{sl}(2)$:
\begin{equation}\label{cocycle}c(X_h,f,g) = \sum_{i=0}^{k+1}\beta_{i} h' f^{(i)}
g^{(k+1-i)}+\sum_{i=0}^{k}\gamma_{i} h'' f^{(i)} g^{(k-i)},\end{equation}
where    $\beta_{i}$ and $\gamma_{i}$ are
constants satisfying:
\begin{equation}\label{coef}
2(\mu-\lambda-\nu-k-1)\gamma_{i}+(i+1)(i+2\lambda)\beta_{i+1}+(k+1-i)(k-i+2\nu)\beta_{i}=0.
\end{equation}
\end{lem}

\begin{proofname}. Any $(k+2)$-homogeneous 1-cocycle on $\frak{sl}(2)$ should retains the following
general form: $$c(X_h,f,g) = \sum_{i=0}^{k+2}\alpha_{i} h
f^{(i)} g^{(k+2-i)} + \sum_{i=0}^{k+1}\beta_{i} h'
f^{(i)} g^{(k+1-i)}+\sum_{i=0}^{k}\gamma_{i} h''
f^{(i)} g^{(k-i)},$$ where $\alpha_{i}$, $\beta_{i}$
and $\gamma_{i}$ are, a priori, functions. First, we prove that
the terms in $h$ can be annihilated by adding a coboundary.
Let $b : \mathcal{F_\lambda}\times \mathcal{F_\nu}\rightarrow
\mathcal{F_\mu}$ be a bilinear differential operator defined by
$$
b(f,g) =
\sum_{i=0}^{k+2}b_{i} f^{(i)} g^{(k+2-i)} ,
$$
where $f\in\mathcal{F_\lambda}$,  $g\in\mathcal{F_\nu}$ and the coefficients $b_{i}$
are functions satisfying $$\begin{array}{l}{d\over d x}(b_{i})=\alpha_{i}.\end{array}$$ Then, for all $X_h
\in\frak{sl}(2)$, we have
\begin{equation*}\begin{array}{lll}\label{coboundary}
\delta b(X_h,f,g) =&
{\displaystyle\sum_{i=0}^{k+2}}\alpha_{i}hf^{(i)}
g^{(k+2-i)}+{\displaystyle\sum_{i=0}^{k+2}}(\mu-\lambda-\nu-k-1)b_{i}h'f^{(i)}
g^{(k+2-i)}\\[10pt]
&-{1\over2}{\displaystyle\sum_{i=0}^{k+1}}\big((i+1)(i+2\lambda)b_{i+1}+(k+2-i)(k+1-i+2\nu)b_{i,}\big)h''f^{(i)}
g^{(k+1-i)}. \end{array}\end{equation*} We replace $c$ by
$\widetilde{c}=c-\delta b$ and then  we see that the
1-cocycle $\widetilde{c}$ does not contain terms in $h$. So, up to a
coboundary, any $(k+2)$-homogeneous 1-cocycle on $\frak{sl}(2)$ can be expressed as
follows:
$$
c(X_h,f,g) = \sum_{i=0}^{k+1}\beta_{i} h'
f^{(i)} g^{(k+1-i)}+\sum_{i=0}^{k}\gamma_{i} h''
f^{(i)} g^{(k-i)}.
$$
Now, consider the 1-cocycle condition:
$$
c([X_{h_1}, X_{h_2}],f,g)- X_{h_1}\cdot c(X_{h_2}, f,g) +X_{h_2}\cdot c(X_{h_1}, f,g) = 0
$$
where $f\in\mathcal{F_\lambda}$,  $g\in\mathcal{F_\nu}$
and  $X_{h_1},\,X_{h_2} \in\frak{sl}(2)$. A direct computation proves that we have
$$
\begin{array}{l}{d\over d x}(\beta_{i})={d\over d x}(\gamma_{m})=0\end{array}
$$
and
\begin{equation*}\label{condition}
2(\mu-\lambda-\nu-k-1)\gamma_{i}+(i+1)(i+2\lambda)\beta_{i+1}+(k+1-i)(k-i+2\nu)\beta_{i}=0.
\end{equation*}
\hfill$\Box$\end{proofname}
\begin{corollary}\label{cor} If $\mu-\l-\nu\neq k+1$, where $k+1\in\mathbb{N}$, then any $(k+2)$-homogeneous 1-cocycle $c\in Z^1_\mathrm{diff}(\frak
{sl}(2),\mathrm{D}_{\lambda,\nu;\mu})$ is a coboundary. Especially,
if $\mu-\l-\nu= k+1$ then any  1-cocycle $c\in
Z^1_\mathrm{diff}(\frak {sl}(2),\mathrm{D}_{\lambda,\nu;\mu})$ is,
up to a coboundary, $(k+2)$-homogeneous and if
$\mu-\l-\nu\notin\mathbb{N}$ then ${\rm H}^1_\mathrm{diff}(\frak
{sl}(2),\mathrm{D}_{\lambda,\nu;\mu})=0$.
\end{corollary}
\begin{proofname}. If $\mu-\l-\nu\neq k+1$ we can easily show that the 1-cocycle $c$ defined by (\ref{cocycle}) is nothing but the operator
$\delta b$ where
$$b(f,g)={1\over \mu-\lambda-\nu-k-1}\sum_{i=0}^{k+1}\beta_{i}
f^{(i)} g^{(k+1-i)}.$$
\hfill$\Box$\end{proofname}

\begin{thm}
\label{th2}(\cite{b})~~
\begin{equation*}
{\rm H}^1_\mathrm{diff}(\frak
{sl}(2),\mathrm{D}_{\lambda,\nu;\mu})\simeq \left\{
\begin{array}{llllll}
\mathbb{R}^3 & \hbox{ if }\quad
 (\lambda,\nu,\mu)\,\text{ is resonant },
\\[2pt]
\mathbb{R} & \hbox{ if }(\lambda,\nu,\mu)\,\text{ is weakly resonant
},
\\[2pt]
0 & \hbox{otherwise}.
\end{array}
\right.
\end{equation*}
\end{thm}

\begin{proofname}.  Let $\mu-\lambda-\nu=k+1$; where $k+1\in\mathbb{N}$, then,
according to Corollary \ref{cor}, any $n$-homogeneous 1-cocycle
$c\in Z^1_\mathrm{diff}(\frak
{sl}(2),\mathrm{D}_{\lambda,\nu;\mu})$, where $n\neq k+2$, is a
coboundary. Thus, we consider only the $(k+2)$-homogeneous
1-cocycles given by Lemma \ref{l33}.  In this case, the relation
(\ref{coef}) becomes:
\begin{equation}\label{beta}
(i+1)(i+2\lambda)\beta_{i+1}+(k+1-i)(k-i+2\nu)\beta_{i}=0.
\end{equation}
Let
$$b(f,g) = \sum_{i=0}^{k+1}b_{i} f^{(i)} g^{(k+1-i)}.
$$ By a direct computation we have
\begin{align*}
\delta b(X_h,f,g)
=-{1\over2}\sum_{i=0}^{k}\big((i+1)(i+2\lambda)b_{i+1}+(k+1-i)(k-i+2\nu)b_{i}\big)h''f^{(i)}
g^{(k-i)}. \end{align*} So, we are in position to complete the proof  as
Bouarroudj did in \cite{b}. We recall here the (slightly modified)
explicit expressions of the 1-cocycles given in \cite{b}. Hereafter,
$\begin{pmatrix}x\\i\end{pmatrix}$ is the standard binomial coefficient:
$\begin{pmatrix}x\\i\end{pmatrix}=\frac{x(x-1)\cdots(x-i+1}{i!}$
that makes sense for arbitrary  $x\in\mathbb{R}$.

{\bf Case 1: $(\lambda,\nu,\mu)$ is weakly resonant}. In this case,
the corresponding cohomology space is one-dimensional, generated by
the 1-cocycle $\mathfrak{a}$ defined as follows:

(i) If $\lambda\neq -{s\over 2}$, where
$s\in\{0,\dots,\, k \}$, then
\begin{equation}\label{a1}
\mathfrak{a}(X_h, f, g) =  \sum_{i=0}^{k+1}
\begin{pmatrix}{k+1}\\i\end{pmatrix}\begin{pmatrix}{2\nu+k}\\i\end{pmatrix}
\begin{pmatrix}{-2\lambda}\\i\end{pmatrix}^{-1}h'f^{(i)}g^{(k+1-i)}.
\end{equation}

(ii)  If $\nu\neq-{t\over 2}$, where
$t \in\{0,\dots,\, k \}$, then
\begin{equation}\label{a2}
\mathfrak{a}(X_h, f, g) =  \sum_{i=0}^{k+1}
\begin{pmatrix}{k+1}\\i\end{pmatrix}\begin{pmatrix}{2\lambda+k}\\k+1-i\end{pmatrix}
\begin{pmatrix}{-2\nu}\\k+1-i\end{pmatrix}^{-1}h'f^{(i)}g^{(k+1-i)}.
\end{equation}

(iii) If $\lambda=-{s\over 2}$ and $\nu= -{t\over 2}$ , where $s,\, t
\in\{0,\dots,\, k \}$ but $s+t< k$, then
\begin{equation}\label{a4}\begin{array}{l}
\mathfrak{a}(X_h, f, g) = {\displaystyle\sum_{i=s+1}^{k-t}}(-1)^{i}\begin{pmatrix}{k+1}\\i\end{pmatrix}
\begin{pmatrix}{k-t-s-1}\\i-s-1\end{pmatrix}h'f^{(i)}g^{(k+1-i)}.
\end{array}
\end{equation}
Observe that if $\mu-\lambda-\nu=0$ then $(\lambda,\nu,\mu)$
is weakly resonant since $\mu-\lambda-\nu\in\mathbb{N}$ but $\mu-\lambda-\nu-1\notin\mathbb{N}$.
In this case, the set $\{0,\dots,\, k \}$ is empty, so we are in  the situations (i) and (ii) and
the 1-cocycle $\mathfrak{a}$ is then defined by
$\mathfrak{a}(X_h,f,g)=h'fg$.\\

{\bf Case 2: $(\lambda,\nu,\mu)$ is  resonant}. That is,  $\lambda=-{s\over 2}$ and $\nu= -{t\over 2}$ , where $s,\, t
\in\{0,\dots,\, k \}$ with $s+t\geq k$. In this case, the
corresponding cohomology space is three-dimensional, generated by
the 1-cocycles $\mathfrak{b}$, $\mathfrak{c}$ and $\mathfrak{d}$
defined as follows:
\begin{equation}\label{b2}\begin{array}{l}
\mathfrak{b}(X_h, f, g) = h''f^{(k-t)} g^{(t)},
\end{array}
\end{equation}
\begin{equation}\label{c2}
\mathfrak{c}(X_h, f, g) = \sum_{i=0}^{s}
\begin{pmatrix}{k+1}\\i\end{pmatrix}
\begin{pmatrix}{k-t}\\i\end{pmatrix}\begin{pmatrix}{s}\\i\end{pmatrix}^{-1}h'f^{(i)}g^{(k+1-i)},
\end{equation}
\begin{equation}\label{d2}
\mathfrak{d}(X_h, f, g) =\sum_{i=s+1}^{k+1}
\begin{pmatrix}{k+1}\\i\end{pmatrix}
\begin{pmatrix}{k-s}\\k+1-i\end{pmatrix}\begin{pmatrix}{t}\\k+1-i\end{pmatrix}^{-1}h'f^{(i)}g^{(k+1-i)}.
\end{equation}
Observe that if $(\lambda,\nu)= (-{s\over 2},-{k-s\over 2})$, where
$s\in\{0,\dots,\, k \}$, then %we  can replace $\mathfrak{c}$ by $\mathfrak{c}+\mathfrak{d}$:
$
(\mathfrak{c}+\mathfrak{d})(X_h, f, g)=h'(fg)^{(k+1)}.$

\hfill$\Box$\end{proofname}

 %%%%%%%%%%%%%%%%%%%%%%%%%%%%%%%%%%%%%%%%%%%%%%%%%%%%%%%%%%%%%%%%%%%%%
 \subsection{Relationship between $\HH^1_\mathrm{diff}(\frak
{osp}(1|2),\frak{D}_{\l,\nu;\mu})$ and $\HH^1_\mathrm{diff}(\frak
{sl}(2),\mathrm{D}_{\lambda,\nu;\mu})$}
 %%%%%%%%%% %%%%%%%%%%%%%%%%%%%%%%%%%%%%%%%%%%%%%%%%%%%%%%%%%%%%%%%%%%%
We need to present here some results illustrating the analogy
between the cohomology spaces in super and classical settings.

\begin{proposition}\label{prop1}~

\noindent 1) As a $\frak {sl}(2)$-module, we have
$$\fF_{\l}\simeq\cF_\l \oplus \Pi(\cF_{\l+\half})\quad\text{ and }\quad \frak
{osp}(1|2)\simeq\frak {sl}(2)\oplus \Pi(\mathcal{H}),$$ where
$\mathcal{H}$ is the subspace of $\cF_{-\half}$ spanned by
$\{dx^{-\frac{1}{2}}, xdx^{-\frac{1}{2}}\}$ and $\Pi$ is the
change of parity.

\noindent 2) As a $\frak {sl}(2)$-module, we have for
the homogeneous relative parity components:
\begin{equation}\label{deven}(\frak{D}_{\lambda,\nu;\mu})_{\bar0}
\simeq\mathrm{D}_{\lambda,\nu;\mu}
\oplus\mathrm{D}_{\lambda+\half,\nu+\half;\mu}
\oplus\mathrm{D}_{\lambda,\nu+\half;\mu+\half}\oplus\mathrm{D}_{\lambda+\half,\nu;\mu+\half},\end{equation}
\begin{equation}\label{dodd}
(\frak{D}_{\lambda,\nu;\mu})_{\bar1}\simeq\Pi\left(\mathrm{D}_{\lambda,\nu;\mu+\half}
\oplus \mathrm{D}_{\lambda+\half,\nu+\half;\mu+\half}
\oplus\mathrm{D}_{\lambda,\nu+\half;\mu}\oplus\mathrm{D}_{\lambda+\half,\nu;\mu}\right).
\end{equation}
\end{proposition}
\begin{proofname}. 1) The first statement is immediately deduced
 from (\ref{deriv}).

2) It is well known that if $M=M_{\bar0}\oplus M_{\bar1}$ and
$N=N_{\bar0}\oplus N_{\bar1}$ are two $\mathfrak{g}$-modules, where
$\mathfrak{g}$ is a (super)algebra, then $\mathrm{Hom}(M,N)$ is a
$\mathfrak{g}$-module, where the homogenous components are
$$
\mathrm{Hom}(M,N)_{\bar0}=\mathrm{Hom}(M_{\bar0},M_{\bar0})
\oplus\mathrm{Hom}(M_{\bar1},N_{\bar1})\quad\text{and}
\quad\mathrm{Hom}(M,N)_{\bar1}=\mathrm{Hom}(M_{\bar0},N_{\bar1})
\oplus\mathrm{Hom}(M_{\bar1},N_{\bar0})
$$
and the $\mathfrak{g}$-action on $\mathrm{Hom}(M,N)$ is given by
\begin{equation*}\label{d-action}
(X.A)(x)=X.(A(x))-(-1)^{|A||X|} A(X.x).
\end{equation*}
Moreover, if $\varphi_1:M\rightarrow M'$ and
$\varphi_2:N\rightarrow N'$ are two $\mathfrak{g}$-isomorphisms,
then the map
$\Psi:\mathrm{Hom}(M,N)\rightarrow\mathrm{Hom}(M',N')$ defined by
\begin{equation*}\label{isop}
\Psi(A)=\varphi_2\circ A\circ\varphi_1^{-1}
\end{equation*}
is a $\mathfrak{g}$-isomorphism.
In our situation, as a $\frak {sl}(2)$-module, we have for the
homogeneous relative parity components:
\begin{equation*}\left\{\begin{array}{llll}
(\frak{F}_{\lambda}\otimes\frak{F}_{\nu})_{\bar0}&\simeq&{\cal
F}_{\lambda}\otimes\mathcal{F}_\nu \oplus\Pi({\cal
F}_{\lambda+\half})\otimes\Pi(\mathcal{F}_{\nu+\half}),
\\(\frak{F}_{\lambda}\otimes\frak{F}_{\nu})_{\bar1}&\simeq&\Pi({\cal
F}_{\lambda+\half})\otimes\mathcal{F}_\nu \oplus{\cal
F}_{\lambda}\otimes\Pi(\mathcal{F}_{\nu+\half}).\end{array}\right.
\end{equation*}
So, we deduce the two homogenous relative parity components of
$\frak D_{\lambda,\nu;\mu}$ as a $\frak {sl}(2)$-module. In fact,
we have the following isomorphisms:
\begin{gather*}\left\{\begin{array}{llllllll}
\mathrm{Hom_{diff}}\left(\Pi(\mathcal{F}_{\lambda+\frac{1}{2}})\otimes
\Pi(\mathcal{F}_{\nu+\frac{1}{2}}),\;\mathcal{F}_\mu\right)
&\rightarrow&\mathrm{D}_{\lambda+\half,\nu+\half;\mu},\quad
&A&\mapsto& A\circ(\Pi\otimes\Pi),\\[4pt]
\mathrm{Hom_{diff}}\left(\mathcal{F}_{\lambda}\otimes
\Pi(\mathcal{F}_{\nu+\frac{1}{2}}),\;\Pi(\mathcal{F}_{\mu+\frac{1}{2}})\right)
&\rightarrow&\mathrm{D}_{\lambda,\nu+\half;\mu+\half},\quad
&A&\mapsto& \Pi\circ A\circ (Id\otimes\Pi),\\[4pt]
\mathrm{Hom_{diff}}\left(\Pi(\mathcal{F}_{\lambda+\frac{1}{2}})\otimes
\mathcal{F}_{\nu},\;\Pi(\mathcal{F}_{\mu+\frac{1}{2}})\right)
&\rightarrow&\mathrm{D}_{\lambda+\half,\nu;\mu+\half}, \quad
&A&\mapsto& \Pi\circ A\circ (\Pi\otimes Id).
\end{array}\right.
\end{gather*}
\begin{gather*}\left\{\small{\begin{array}{llllllll}
\mathrm{Hom_{diff}}\left(\mathcal{F}_{\lambda}\otimes
\mathcal{F}_{\nu},\Pi(\mathcal{F}_{\mu+\frac{1}{2}})\right)
&\rightarrow\Pi(\mathrm{D}_{\lambda,\nu;\mu+\half}), &A\mapsto& \Pi(\Pi\circ A),\\[4pt]
\mathrm{Hom_{diff}}\left(\Pi(\mathcal{F}_{\lambda+\half})\otimes
\Pi(\mathcal{F}_{\nu+\frac{1}{2}}),\Pi(\mathcal{F}_{\mu+\frac{1}{2}})\right)
&\rightarrow\Pi(\mathrm{D}_{\lambda+\half,\nu+\half;\mu+\half}),
&A\mapsto& \Pi(\Pi\circ A\circ (\Pi\otimes\Pi)),\\[4pt]
\mathrm{Hom_{diff}}\left(\mathcal{F}_{\lambda}\otimes
\Pi(\mathcal{F}_{\nu+\frac{1}{2}}),\mathcal{F}_{\mu}\right)
&\rightarrow\Pi(\mathrm{D}_{\lambda,\nu+\half;\mu}),
&A\mapsto& \Pi(A\circ (Id\otimes \Pi)),\\[4pt]
\mathrm{Hom_{diff}}\left(\Pi(\mathcal{F}_{\lambda+\frac{1}{2}})\otimes
\mathcal{F}_{\nu},\mathcal{F}_{\mu}\right) &\rightarrow\Pi({
D}_{\lambda+\half,\nu;\mu}), &A\mapsto& \Pi(A\circ (\Pi\otimes Id)).
\end{array}}\right.
\end{gather*}
\hfill$\Box$\end{proofname}

Now, in order to compute $\HH^1(\frak
{osp}(1|2),\frak{D}_{\lambda,\nu;\mu})$, we need first to describe the $\mathfrak{sl}(2)$-relative cohomology
space $\HH^1_\mathrm{diff}(\mathfrak{osp}(1|2),\frak
{sl}(2);\mathfrak{D}_{\lambda,\nu;\mu})$. So, we shall need the following
description of some $\mathfrak{sl}(2)$-invariant
 mappings.
\begin{lem}
\label{inva}Let
\begin{equation*}A:\mathcal{H}\otimes\cF_\l\otimes\cF_\nu\rightarrow\cF_\mu,\quad
(h(dx)^{-\half},f(dx)^\l,g(dx)^\nu)\mapsto
A(h,f,g)(dx)^{\mu}\end{equation*} be a trilinear differential
operator. If $A$ is a nontrivial  $\frak{sl}(2)$-invariant operator then
\begin{equation*}\begin{array}{ll}\mu=\l+\nu+k-\half,\quad\text{where}\quad
k\in\mathbb{N}.\end{array}\end{equation*}  For $k\geq2$, the corresponding
operator $A_k$ is  of the form:
\begin{equation*}
A_k(h,f,g)=\sum_{i=0}^kc_{i}hf^{(i)}g^{(k-i)}+\sum_{i=0}^{k-1}\big[(i+1)(i+2\lambda)c_{i+1}+
(k-i)(k-i-1+2\nu)c_{i})\big]h'f^{(i)}g^{(k-i-1)},
\end{equation*}where the $c_{i}$ are constant characterized by
the following recurrence formula:
\begin{equation}\small{\label{rec}\begin{array}{ll}
&(i+1)(i+2)(i+2\lambda)(i+2\lambda+1)c_{i+2}+2(i+1)(k-i-1)(i+2\lambda)(k-i-2+2\nu)c_{i+1}\\[4pt]&~~~~~~~~~~~~+
(k-i-1)(k-i)(k-i-2+2\nu)(k-i-1+2\nu)c_{i}=0.\end{array}}
\end{equation}
For $k=0,\,1$, we have
$$A_0(h,f,g)=c_{0}hfg\quad\text{and}\quad A_1(h,f,g)=c_{0}hfg+c_1hf'g+(2\lambda c_1+2\nu c_0)h'fg.$$
\end{lem}
\begin{proofname}. Obviously, the operator $A$ is $\frak{sl}(2)$-invariant if and only if each of its homogenous components is $\frak{sl}(2)$-invariant. Moreover, the invariance with respect the vector field $X_1=\partial_x$ yields that $A$ must be expressed with constant coefficients. Thus, let $k\in\mathbb{N}$ and consider \begin{equation*}
A_k(h,f,g)=\sum_{i=0}^kc_{i}hf^{(i)}g^{(k-i)}+\sum_{i=0}^{k-1}d_{i}h'f^{(i)}g^{(k-i-1)},
\end{equation*}where the $c_{i}$ and $d_{i}$ are constants. The
invariance property of $A$ with respect any vector fields $X_F$ reads:
\begin{equation}\label{inv}\begin{array}{llll} F(A(h,f,g))'+\mu F'A(h,f,g)&=&A(Fh'-{1\over2}F'h,f,g)+
A(h,Ff'+\lambda F'f,g)\\&~&~~~~~~~~~~~~~~~~~~~~~~~~~~
+A(h,f,Fg'+\nu F'g).\end{array}\end{equation}
Consider any non vanishing coefficient $c_{i}$ and consider terms in
$F'hf^{(i)}g^{(k-i)}$ in (\ref{inv}), we get
\begin{equation*}\begin{array}{ll}\mu=\l+\nu+k-\half.\end{array}\end{equation*}
Considering respectively terms in $F''hf^{(i)}g^{(k-i-1)}$ and (for
$k\geq2$) $F''h'f^{(i)}g^{(k-i-2)}$ yield
\begin{equation}\label{cij}d_{i}=(i+1)(i+2\lambda)c_{i+1}+ (k-i)(k-i-1+2\nu)c_{i}\end{equation}
\begin{equation}\label{dij}0=(i+1)(i+2\lambda)d_{i+1}+ (k-i)(k-i-1+2\nu)d_{i}.\end{equation}
Combining (\ref{cij}) and (\ref{dij}) we have (\ref{rec}). Under
these conditions we check that the operator $A_k$ is
$\frak{sl}(2)$-invariant.
\hfill$\Box$ \end{proofname}

\begin{prop}\label{sa} The $\mathfrak{sl}(2)$-relative cohomology
spaces $\HH^1_\mathrm{diff}(\mathfrak{osp}(1|2),\frak
{sl}(2);\mathfrak{D}_{\lambda,\nu;\mu})$ are all trivial. That is,
any 1-cocycle $\Upsilon$ is a coboundary over $\frak{osp}(1|2)$ if
and only if its restriction  to $\frak{sl}(2)$ is a coboundary over
$\frak{sl}(2)$.
\end{prop}
\begin{proofname}. First, it is easy to see that any 1-cocycle $\Upsilon\in
Z^1_{\mathrm{diff}}(\mathfrak{osp}(1|2),\frak{D}_{\lambda,\nu;\mu})$
vanishing on $\frak{sl}(2)$ is $\mathfrak{sl}(2)$-invariant. Indeed,
the 1-cocycle relation of $\Upsilon$ reads:
\begin{equation*}\label{osp1}
(-1)^{|F||\Upsilon|}X_F\cdot
\Upsilon(X_G)-(-1)^{|G|(|F|+|\Upsilon|)}X_G\cdot
\Upsilon(X_F)-\Upsilon([X_F,~X_G])=0,
\end{equation*}
where $X_F,\,X_G\in \mathfrak{osp}(1|2).$  If $\Upsilon(X_F)=0$ for
all $X_F\in\frak{sl}(2)$ then the previous equation  becomes
\begin{equation*}\label{osp2}
X_F\cdot \Upsilon(X_G)-\Upsilon([X_F,~X_G])=0
\end{equation*}
 expressing the $\frak{sl}(2)$-invariance  of
$\Upsilon$. Thus, the space
$\HH^1_\mathrm{diff}(\mathfrak{osp}(1|2),\frak
{sl}(2);\mathfrak{D}_{\lambda,\nu;\mu})$ is nothing but the space of
cohomology classes of 1-cocycles vanishing on $\frak{sl}(2)$.

Let $\Upsilon$ be a 1-cocycle vanishing on $\frak{sl}(2)$, then, by the 1-cocycle condition, we have:
\begin{align}
\label{sltr1} &{X_f}\cdot\Upsilon(X_{h\theta})-\Upsilon([X_f,X_{h\theta}])=0, \\[2pt]
\label{sltr2}&{X_{h_1\theta}}\cdot\Upsilon(X_{h_2\theta})+
{X_{h_2\theta}}\cdot\Upsilon(X_{h_1\theta})=0,
 \end{align}
where  $f\in\mathbb{R}_2[x]$  and $h$, $h_1$,
$h_2\in\mathbb{R}_1[x]$. Here, $\mathbb{R}_n[x]$ is the space of
polynomial functions in the variable $x$, with degree at most $n$.

1) If $\Upsilon$ is an even 1-cocycle then $\Upsilon$ is decomposed
into four trilinear maps:
\begin{gather*}\left\{\begin{array}{llllllll}
\Pi(\mathcal{H})\otimes\Pi(\mathcal{F}_{\lambda+\frac{1}{2}})\otimes
\mathcal{F}_{\nu}
&\rightarrow&{\cal F}_{\mu},\\[4pt]
\Pi(\mathcal{H})\otimes\mathcal{F}_{\lambda}\otimes
\Pi(\mathcal{F}_{\nu+\frac{1}{2}})
&\rightarrow&{\cal F}_{\mu},\\[4pt]
\Pi(\mathcal{H})\otimes\mathcal{F}_{\lambda}\otimes
\mathcal{F}_{\nu}
&\rightarrow&\Pi({\cal F}_{\mu+\frac{1}{2}}),\\[4pt]
\Pi(\mathcal{H})\otimes\Pi(\mathcal{F}_{\lambda+\frac{1}{2}})\otimes
\Pi(\mathcal{F}_{\nu+\frac{1}{2}}) &\rightarrow&\Pi({\cal
F}_{\mu+\frac{1}{2}}).
\end{array}\right.
\end{gather*}
The equation (\ref{sltr1}) is nothing but the
$\frak{sl}(2)$-invariance property of these maps. Therefore, the
expressions of these maps are given by Lemma \ref{inva}, in fact,
the change of parity functor $\Pi$ commutes with the the
$\mathfrak{sl}(2)$-action. So,  we must have $\mu=\l+\nu+k,$ where
$k+1\in\mathbb{N}$, otherwise, the operator $\Upsilon$ is
identically the zero map. More precisely:

If $\mu=\l+\nu+k$ where $k\in\mathbb{N}^*$, we have
\begin{equation}\label{eq1} \Upsilon_k(X_{\theta h})(\theta
f,g)=\sum_{i=0}^{k}c^1_{i}hf^{(i)}g^{(k-i)}+
\sum_{i=0}^{k-1}d^1_{i}h'f^{(i)}g^{(k-i-1)},\end{equation}
\begin{equation}
\Upsilon_k(X_{\theta h})(f,\theta
g)=\sum_{i=0}^kc^2_{i}hf^{(i)}g^{(k-i)}+
\sum_{i=0}^{k-1}d^2_{i}h'f^{(i)}g^{(k-i-1)},\end{equation}
\begin{equation}
\Upsilon_k(X_{\theta
h})(f,g)=\theta\sum_{i=0}^{k+1}c^3_{i}hf^{(i)}g^{(k-i+1)}+
\theta\sum_{i=0}^{k}d^3_{i}h'f^{(i)}g^{(k-i)},\end{equation}
\begin{equation}\label{eq4}
\Upsilon_k(X_{\theta h})(\theta f,\theta
g)=\theta\sum_{i=0}^{k}c^4_{i}hf^{(i)}g^{(k-i)}+
\theta\sum_{i=0}^{k-1}d^4_{i}h'f^{(i)}g^{(k-i-1)},\end{equation}
where
\begin{align*}d^1_{i}&=(i+1)(i+2\lambda+1)c^1_{i+1}+
(k-i)(k-i-1+2\nu)c^1_{i},\\
d^2_{i}&=(i+1)(i+2\lambda)c^2_{i+1}+
(k-i)(k-i-1+2\nu+1)c^2_{i,},\\
d^3_{i}&=(i+1)(i+2\lambda)c^3_{i+1}+
(k-i)(k-i-1+2\nu)c^3_{i},\\
d^4_{i}&=(i+1)(i+2\lambda+1)c^4_{i+1}+
(k-i)(k-i-1+2\nu+1)c^4_{i}\end{align*} and where the coefficients $c^r_{i}$ are satisfying
the recurrence formulae (\ref{rec}).

If $\mu=\l+\nu-1$, we have
\begin{equation*}{\Upsilon}_{-1}(X_{\theta h})(\theta f,\,g)={\Upsilon}_{-1}(X_{\theta h})(f,\,\theta
g)={\Upsilon}_{-1}(X_{\theta h})(\theta f,\,\theta
g)=0,\end{equation*}
\begin{equation*}{\Upsilon}_{-1}({X_\theta h})(f,\,g)=c^{3}_{0}\theta hfg.\end{equation*}

If $\mu=\l+\nu$, we have
\begin{equation*}{\Upsilon}_0(X_{\theta h})(\theta f,\,g)=c^{1}_{0}hfg,\quad{\Upsilon}_0(X_{\theta h})(f,\,\theta
g)=c^{2}_{0}hfg,\quad{\Upsilon}_0(X_{\theta h})(\theta f,\,\theta
g)=\theta c^{4}_{0}hfg,\end{equation*}
\begin{equation*}{\Upsilon}_0(X_{\theta h})(f,\,g)=\theta[c^{3}_{0}hfg'+c^{3}_{1}hf'g+(2\lambda
c^{3}_{1}+2\nu c^{3}_{0})h'fg].\end{equation*}

The maps $\Upsilon_k$ must satisfy the equation (\ref{sltr2}). More
precisely, the maps $\Upsilon_k$ satisfy the following four
equations
\begin{equation*}\label{tfg}\begin{array}{ll}
  {1\over2}\theta
h_1(\Upsilon_k(X_{h_2\theta})(f\theta,g))'&+\mu\theta
h'_1\Upsilon_k(X_{h_2\theta})(f\theta,g)+
\Upsilon_k(X_{h_2\theta})({1\over2}h_1f,g)\\[10pt]&-\Upsilon_k(X_{h_2\theta})(f\theta,\theta({1\over2}h_1g'+\nu
h'_1g))+(h_1\leftrightarrow h_2)=0,
 \end{array}\end{equation*}
 \begin{equation*}\label{ftg}\begin{array}{ll}
{1\over2}\theta
h_1(\Upsilon_k(X_{h_2\theta})(f,g\theta))'&+\mu\theta
h'_1\Upsilon_k(X_{h_2\theta})(f,g\theta)+\Upsilon_k(X_{h_2\theta})(\theta({1\over2}h_1f'+\lambda
h'_1f),g\theta)\\[10pt]&+\Upsilon_k(X_{h_2\theta})(f,{1\over2}h_1g)+(h_1\leftrightarrow
h_2)=0,
 \end{array}\end{equation*}
\begin{equation*}\label{fg}\begin{array}{ll}{1\over2}h_1\partial_\theta
(\Upsilon_k(X_{h_2\theta})(f,g))&+\Upsilon_k(X_{h_2\theta})(\theta({1\over2}h_1f'+\lambda
h'_1f) ,g)\\[10pt]&+\Upsilon_k(X_{h_2\theta})(f,\theta({1\over2}h_1g'+\nu
h'_1g))+(h_1\leftrightarrow h_2)=0,
 \end{array}\end{equation*}
 \begin{equation*}\label{tftg}\begin{array}{ll}
{1\over2}h_1\partial_\theta
(\Upsilon_k(X_{h_2\theta})(f\theta,g\theta))&+\Upsilon_k(X_{h_2\theta})({1\over2}h_1f,g\theta)
\\[10pt]&-\Upsilon_k(X_{h_2\theta})(f\theta,{1\over2}h_1g))+(h_1\leftrightarrow
h_2)=0.
 \end{array}\end{equation*}

 By a direct, but very hard, computation we show that $\Upsilon_k$ is a coboundary. For instance, if $\nu,
\lambda\not\in\{0,\,-\half,\,-1,\,\dots,\,-{k\over2}\}$, we check that $\Upsilon_k=\delta B_k$
where
$$
B_k(f_0+f_1\theta,\,g_0+g_1\theta)=\theta\sum_i(-1)^{i}\begin{pmatrix}{k-1}\\i\end{pmatrix}
\begin{pmatrix}2\nu+k-1\\i\end{pmatrix}
\begin{pmatrix}2\lambda+i\\i\end{pmatrix}^{-1}f_1^{(i)}g_0^{(k-i)}.
$$
\vskip .3cm 2) Similarly, if $\Upsilon$ is an odd 1-cocycle then
$\Upsilon$ is decomposed into four components:
\begin{gather*}\left\{\begin{array}{llllllll}
\Pi(\mathcal{H})\otimes{\cal F}_{\lambda}\otimes
\mathcal{F}_{\nu}
&\rightarrow&{\cal F}_{\mu},\\[4pt]
\Pi(\mathcal{H})\otimes\Pi(\mathcal{F}_{\lambda+\frac{1}{2}})\otimes
\Pi(\mathcal{F}_{\nu+\frac{1}{2}})
&\rightarrow&{\cal F}_{\mu},\\[4pt]
\Pi(\mathcal{H})\otimes\Pi(\mathcal{F}_{\lambda+\frac{1}{2}})\otimes
\mathcal{F}_{\nu}
&\rightarrow&\Pi({\cal F}_{\mu+\frac{1}{2}}),\\[4pt]
\Pi(\mathcal{H})\otimes{\cal F}_{\lambda}\otimes
\Pi(\mathcal{F}_{\nu+\frac{1}{2}}) &\rightarrow&\Pi({\cal
F}_{\mu+\frac{1}{2}}).
\end{array}\right.
\end{gather*}
The equation (\ref{sltr1}) is nothing but the
$\frak{sl}(2)$-invariance property of these bilinear maps. Therefore, the
expressions of these maps are given by Lemma \ref{inva}. So,  we must have
$\mu=\l+\nu+k-\frac{1}{2},$ where $k\in\mathbb{N}$, otherwise, the
operator $\Upsilon$ is identically the zero map. If
$\mu=\l+\nu+k-\frac{1}{2}$ where $k\in\mathbb{N}$, we show, as in the previous case that $\Upsilon$ is a coboundary.

\hfill$\Box$
\end{proofname}
\begin{lemma}\label{cor1}
\label{sd}  Up to a coboundary, any 1-cocycle $\Upsilon\in Z^1_\mathrm{diff}(\frak{
{osp}}(1|2),\frak{D}_{\lambda,\nu;\mu})$ is invariant with respect the vector field $X_1={\partial_x}$. That is, the map $\Upsilon$ can be expressed with constant coefficients.
\end{lemma}

\begin{proofname}. The 1-cocycle condition reads:
\begin{equation}
\begin{array}{lll}\label{part1}
X_1\cdot\Upsilon(X_F)-
(-1)^{|F||\Upsilon|}{X_F}\cdot\Upsilon(X_1)-
\Upsilon([X_1,X_F])=0.
\end{array}
\end{equation}
But, from  Theorem \ref{th1}, up to a coboundary, we have
$\Upsilon(X_1)=0$, and therefore the equation (\ref{part1}) becomes
\begin{equation*}
\begin{array}{lll}\label{}
X_1\cdot(\Upsilon(X_F))-
\Upsilon([X_1,X_F])=0
\end{array}
\end{equation*}
which is nothing but the invariance property of $\Upsilon$ with
respect the vector field $X_1={\partial_x}$.
\hfill$\Box$\end{proofname}

%%%%%%%%%%%%%%%%%%%%%%%%%%%%%%%%%%%%%%%%%%%%%%%%%%%%%%%%%%%%%%%%%%%%%%%%%%%%%%%%%%%%
\section{Proof of Theorem \ref{th1}}\label{s4}
%%%%%%%%%%%%%%%%%%%%%%%%%%%%%%%%%%%%%%%%%%%%%%%%%%%%%%%%%%%%%%%%%%%%%%%%%%%%%%%%%%%

First, according to Proposition \ref{prop1}, Proposition \ref{sa}
and Theorem \ref{th2}, we easily check that the following statements
hold:
\begin{itemize}
  \item [i)] The space $\HH^1_\mathrm{diff}(\frak
{osp}(1|2),\frak{D}_{\l,\nu;\mu})$ is trivial if
$2(\mu-\l-\nu)+1\notin\mathbb{N}$.
  \item [ii)] The space $\HH^1_\mathrm{diff}(\frak
{osp}(1|2),\frak{D}_{\l,\nu;\mu})$ is even if $\mu-\lambda-\nu$ is integer  and  it is  odd
if $\mu-\lambda-\nu$ is semi-integer.
\end{itemize}

\begin{prop}\label{prop2} Let $\Upsilon\in Z^1_\mathrm{diff}(\frak{
{osp}}(1|2),\frak{D}_{\lambda,\nu;\mu})$,  $k\in\mathbb{N}$, $h\in\mathbb{R}_1[x]$ and $f,\,g\in C^\infty(\mathbb{R})$.

a) If $\mu-\l-\nu=k$ then, up to a coboundary, $\Upsilon(X_{\theta h})(\theta f,\theta
g)$, $\Upsilon(X_{\theta h})( f,\theta
g)$ and $\Upsilon(X_{\theta h})(\theta f,
g)$ are $k$-homogeneous and $\Upsilon(X_{\theta h})(f,
g)$ is $(k+1)$-homogeneous.

b) If $\mu-\l-\nu=k-{1\over2}$ then, up to a coboundary, $\Upsilon(X_{\theta h})(f,
g)$, $\Upsilon(X_{\theta h})( f,\theta
g)$ and $\Upsilon(X_{\theta h})(\theta f,
g)$ are $k$-homogeneous and $\Upsilon(X_{\theta h})(\theta f,\theta
g)$ is $(k-1)$-homogeneous.
\end{prop}
\begin{proofname}. Let $\mu-\l-\nu=k$. Up to a coboundary, the operator
$\Upsilon(X_{\theta h})$  is an odd map. Therefore, the elements
$\Upsilon(X_{\theta h})(\theta f,\theta
g)$, $\Upsilon(X_{\theta h})( f,\theta
g)$, \dots are all homogeneous (even or odd). Thus, the actions of
$X_f$ and $X_{h\theta}$ on these elements are also homogeneous, see
(\ref{actiondecomp}).

Now, for $h=x$ and $f=x^2$, the equation (\ref{sltr1}) becomes
$$
X_{x^2}\cdot
   \Upsilon(X_{x\theta}) =
 X_{x\theta} \cdot\Upsilon(X_{x^2}).
$$
So, using Lemma \ref{cor} and formulas (\ref{actiondecomp}), we
obtain the statement a) for $h=x$. Besides, we use again the
equation (\ref{sltr1}) but for $h=1$ and $f=x$. The statement b) can
be proved similarly. \hfill$\Box$\end{proofname}\\

Now, we explain our strategy to prove Theorem \ref{th1}. Consider
$\Upsilon\in Z^1_\mathrm{diff}(\frak{
{osp}}(1|2),\frak{D}_{\lambda,\nu;\mu})$ where  $2(\mu-\l-\nu)+1\in\mathbb{N}$.
That is,
$$
\begin{array}{l}\mu-\l-\nu=k\quad\text{or}\quad \mu-\l-\nu=k-{1\over2}
\quad\text{where}\quad k\in\mathbb{N}.\end{array}
$$
For instance, in the first case, the cohomology space is even, therefore,
the restriction of $\Upsilon$ on $\frak
{sl}(2)$ is with values in $(\frak{D}_{\lambda,\nu;\mu})_{\bar0}$
which is isomorphic, as $\frak
{sl}(2)$- module,  to
$$
\mathrm{D}_{\lambda,\nu;\mu}
\oplus\mathrm{D}_{\lambda+\half,\nu+\half;\mu} \oplus
\mathrm{D}_{\lambda,\nu+\half;\mu+\half}\oplus
\mathrm{D}_{\lambda+\half,\nu;\mu+\half},
$$
while the restriction of $\Upsilon$ on $\Pi
(\mathcal{H})$ is with values in $(\frak{D}_{\lambda,\nu;\mu})_{\bar1}$  which is isomorphic, as $\frak
{sl}(2)$- module, to
$$
\Pi( \mathrm{D}_{\lambda,\nu;\mu+\half} \oplus
\mathrm{D}_{\lambda+\half,\nu+\half;\mu+\half} \oplus
\mathrm{D}_{\lambda,\nu+\half;\mu}\oplus\mathrm{D}_{\lambda+\half,\nu;\mu}).
$$
Now, according to the decompositions (\ref{deven}) and (\ref{dodd}),
we have
\begin{equation*}\begin{array}{lll}\HH^1(\frak
{sl}(2),\frak{D}_{\l,\nu;\mu})&=\HH^1(\frak {sl}(2),{\rm
D}_{\lambda,\nu;\mu} \oplus{\rm D}_{\lambda+\half,\nu+\half;\mu}
\oplus{\rm D}_{\lambda,\nu+\half;\mu+\half}\oplus{\rm
D}_{\lambda+\half,\nu;\mu+\half})\\[8pt]&~~\oplus\HH^1(\frak
{sl}(2),\Pi({\rm D}_{\lambda,\nu;\mu+\half} \oplus{\rm
D}_{\lambda+\half,\nu+\half;\mu+\half} \oplus{\rm
D}_{\lambda,\nu+\half;\mu}\oplus{\rm
D}_{\lambda+\half,\nu;\mu})).\end{array}\end{equation*} Thus, the
restriction of $\Upsilon$ to $\frak {sl}(2)$ is, a priori, described
by Theorem \ref{th1} while the general form of the restriction of
$\Upsilon$ on $\Pi (\mathcal{H})$ is given by Proposition
\ref{prop2}. Finally, the operator $\Upsilon$ will be completely
given by the 1-cocycle conditions.

Hereafter, $F=f_0+f_1\theta$ and $G=g_0+g_1\theta$ where
$f_0,\,g_0,\,f_1,\,g_1\in C^\infty(\mathbb{R})$.\\

\noindent{\bf Case 1: $(\lambda,\nu,\mu)$ is weakly super resonant}
with $\mu-\lambda-\nu=k+1\in\mathbb{N}$.
In this case, we describe the restriction of
$\Upsilon$ to $\mathfrak{sl}(2)$ by using the 1-cocycles (\ref{a1}), (\ref{a2})
and (\ref{a4}).

a)  Let $(\lambda,\nu)\neq(-\frac{s}{2},-\frac{t}{2})$ where
$s,\,t\in\{0,1,\ldots,k+1\}$. In this case, the 1-cocycle
$\Upsilon$ is even and (if, for instance,
$\lambda\neq-\frac{s}{2}$) its restriction to $\frak {sl}(2)$ is
given by
\begin{equation*}
\Upsilon(X_h,F,G)=\alpha_1\mathfrak{a}_1(X_h,f_0,g_0)+
\alpha_2\mathfrak{a}_2(X_h,f_1,g_1)+\theta\alpha_3\mathfrak{a}_3(X_h,f_1,g_0)+
\theta\alpha_4\mathfrak{a}_4(X_h,f_0,g_1)
\end{equation*}
where
$\alpha_1,\,\alpha_2,\,\alpha_3,\,\alpha_4\in\mathbb{R}$ and the
maps
$\mathfrak{a}_1,\,\mathfrak{a}_2,\,\mathfrak{a}_3,\,\mathfrak{a}_4$
are as in (\ref{a1}). For instance, the expression of
$\mathfrak{a}_2$ can be deduced from (\ref{a1}) by substituting
respectively $\lambda+{1\over2}$, $\nu+{1\over2}$ and $k-1$ to
$\lambda$, $\nu$ and $k$, see (\ref{deven}). From the relation
$\delta\Upsilon(X_{h_0},X_{h_1\theta})(F,G)=0$ we deduce that
$$\begin{array}{l}
\alpha_{4}=\alpha_{1},\quad
2\lambda\,\alpha_{3}=(2\lambda+k+1)\alpha_{1}\quad\text{and}\quad
2\lambda\,\alpha_2=-(k+1)\alpha_{1}.\end{array}
$$
Thus, according to Proposition \ref{sa}, we have
$\mathrm{dim}\mathrm{H}^1(\frak{osp}
({1|2}),(\frak{D}_{\lambda,\nu;\mu}))\leq1$. Now, using Lemma
\ref{cor1}, Proposition \ref{prop2}, the isomorphism (\ref{deven})
and the 1-cocycle relations, we extend $\Upsilon$ to
$\Pi(\mathcal{H})$. More precisely, we prove that we can choose
\begin{equation*}\begin{array}{llll} &\Upsilon(X_{h\theta}, F,G)=\alpha_1\theta h'
\Big({\displaystyle\sum_{i=0}^{k+1}}\begin{pmatrix}{k+1}\\i\end{pmatrix}
\begin{pmatrix}{2\nu+k}\\i\end{pmatrix}\begin{pmatrix}{-2\lambda}\\i\end{pmatrix}^{-1}
f_0^{(i)}g_0^{(k+1-i)}\\[4pt]
&~~~~~~~~~-\frac{(k+1)}{2\lambda}f_1g_1^{(k)}-
{\displaystyle\sum_{i=1}^{k}}\begin{pmatrix}{k}\\i\end{pmatrix}
\begin{pmatrix}{2\nu+k-1}\\i\end{pmatrix}\begin{pmatrix}{-2\lambda}\\i\end{pmatrix}^{-1}
f_1^{(i)}g_1^{(k-i)}\Big).
\end{array}\end{equation*}
Thus, in this case,
$\mathrm{dim}\mathrm{H}^1(\frak{osp}
({1|2}),(\frak{D}_{\lambda,\nu;\mu}))=1$. \\

b) Let
$(\lambda,\nu)=(-\frac{s}{2},-\frac{t}{2})$ with
$s,\,t~\in\{0,\,1,\ldots,k+1\}$ but  $s+t< k$. As in in the previous case, the  restriction
 of $\Upsilon$ to $\frak {sl}(2)$ is given by
\begin{equation*}
\Upsilon(X_h,F,G)=\alpha_1\mathfrak{a}_1(X_h,f_0,g_0)+
\alpha_2\mathfrak{a}_2(X_h,f_1,g_1)+\theta\alpha_3\mathfrak{a}_3(X_h,f_1,g_0)+
\theta\alpha_4\mathfrak{a}_4(X_h,f_0,g_1)
\end{equation*}
where
$\alpha_1,\,\alpha_2,\,\alpha_3,\,\alpha_4\in\mathbb{R}$, but here the
maps
$\mathfrak{a}_1,\,\mathfrak{a}_2,\,\mathfrak{a}_3,\,\mathfrak{a}_4$
are as in (\ref{a4}). By using again the 1-cocycle relations we prove that
$$\begin{array}{l}
\alpha_2=-\frac{k+1}{k-s+1}\,\alpha_3\quad
\alpha_4=-\frac{k-t+1}{k-s+1}\,\alpha_3\quad
\alpha_1=-\frac{k-t-s}{k-s+1}\,\alpha_3.
\end{array}$$
We prove that $\Upsilon$ can be extended  to
$\Pi(\mathcal{H})$. For instance, we can choose
\begin{equation*}\begin{array}{llll} &\Upsilon(X_{h\theta}, F,G)=-\frac{\alpha_3}{k-s+1}\theta h'
\Big({(k-t+1)\displaystyle\sum_{i=s+1}^{k-t}}(-1)^{i}\begin{pmatrix}{k+1}\\i\end{pmatrix}
\begin{pmatrix}{k-t-s-1}\\i-s-1\end{pmatrix}
f_0^{(i)}g_0^{(k-i)}\\[4pt]
&~~~~~~~~~~~~~~~~~~~+
{(k+1)\displaystyle\sum_{i=s}^{k-t}}(-1)^{i}\begin{pmatrix}{k+1}\\i\end{pmatrix}
\begin{pmatrix}{k-t-s-1}\\i-s-1\end{pmatrix}
f_1^{(i)}g_1^{(k-i)}\Big).
\end{array}\end{equation*}
Thus, in this case, $\mathrm{dim}\mathrm{H}^1(\frak{osp}
({1|2}),(\frak{D}_{\lambda,\nu;\mu}))=1$.\\

The case $\mu-\lambda-\nu=k+{3\over2}$ where
$(\lambda,\nu)\neq(-\frac{s}{2},-\frac{t}{2})$ or
$(\lambda,\nu)=(-\frac{s}{2},-\frac{t}{2})$ but $s+t< k+1$ with
$s,\,t~\in\{0,\,1,\ldots,k\}$ can be treated similarly. For
instance, let $(\lambda,\nu)\neq(-\frac{s}{2},-\frac{t}{2})$ where
$s,\,t\in\{0,1,\ldots,k+1\}$. The 1-cocycle $\Upsilon$ is odd and
(if, for instance, $\lambda\neq-\frac{s}{2}$) its restriction to
$\frak {sl}(2)$ is given by
\begin{equation*}
\Upsilon(X_h,F,G)=\theta\alpha_1\mathfrak{a}_1(X_h,f_0,g_0)+
\theta\alpha_2\mathfrak{a}_2(X_h,f_1,g_1)+\alpha_3\mathfrak{a}_3(X_h,f_1,g_0)+
\alpha_4\mathfrak{a}_4(X_h,f_0,g_1)
\end{equation*}
where
$\alpha_1,\,\alpha_2,\,\alpha_3,\,\alpha_4\in\mathbb{R}$ and the
maps
$\mathfrak{a}_1,\,\mathfrak{a}_2,\,\mathfrak{a}_3,\,\mathfrak{a}_4$
are as in (\ref{a1}). For instance, the expression of
$\mathfrak{a}_1$ can be deduced from (\ref{a1}) by substituting $k+1$ to $k$ while the expression of
$\mathfrak{a}_2$ can be deduced from (\ref{a1}) by substituting
respectively $\lambda+{1\over2}$ and $\nu+{1\over2}$ to
$\lambda$ and $\nu$, see (\ref{deven}).
From the relation
$\delta\Upsilon(X_{h_0},X_{h_1\theta})(F,G)=0$ we deduce that
$$\begin{array}{l}
\alpha_{2}=-\frac{2\lambda}{2\nu+k+1}\,\alpha_1, \quad
\alpha_3=\frac{2\nu+2\lambda+k+1}{2\nu+k+1}\alpha_1\quad
\text{and}\quad\alpha_4=\frac{2\lambda}{2\nu+k+1}\,\alpha_1\end{array}
$$
and we prove also that $\Upsilon$ can be extended to $\Pi(\mathcal{H})$.\\

\noindent{\bf Case 2: $(\lambda,\nu,\mu)$ is super resonant}:
$\mu-\lambda-\nu=k+1$ where
$(\lambda,\nu)=(-\frac{s}{2},-\frac{t}{2})$ with
$s,\,t~\in\{1,\ldots,k\}$ and $s+t\geq k+1$. In this case the map
$\Upsilon|_{\mathfrak{sl}(2)}$ can be decomposed as follows:
$\Upsilon|_{\mathfrak{sl}(2)}=B+C+D$ where
\begin{equation*}
B(X_h,F,G)=h''(\beta_1f_0^{(k-t)} g_0^{(t)}+ \beta_2f_1^{(k-t)}
g_1^{(t-1)}+\theta\beta_3f_1^{(k-t)} g_0^{(t)}+
\theta\beta_4f_0^{(k-t+1)} g_1^{(t-1)}),
\end{equation*}\begin{equation*}
C(X_h,F,G)=\gamma_1\mathfrak{c}_1(X_h,f_0,g_0)+
\gamma_2\mathfrak{c}_2(X_h,f_1,g_1)+\theta\gamma_3\mathfrak{c}_3(X_h,f_1,g_0)+
\theta\gamma_4\mathfrak{c}_4(X_h,f_0,g_1),
\end{equation*}\begin{equation*}
D(X_h,F,G)=\delta_1\mathfrak{d}_1(X_h,f_0,g_0)+
\delta_2\mathfrak{d}_2(X_h,f_1,g_1)+\theta\delta_3\mathfrak{d}_3(X_h,f_1,g_0)+
\theta\delta_4\mathfrak{d}_4(X_h,f_0,g_1)
\end{equation*}where the $\mathfrak{c}_i$ and the $\mathfrak{d}_i$
are as those defined in (\ref{c2}) and (\ref{d2}). By the 1-cocycle
relation: $\delta\Upsilon(X_{h_0},X_{h_1\theta})(F,G)=0$ we prove
that
$$\begin{array}{l}
\gamma_{4}=\gamma_{1},\quad \delta_1=\delta_3,\quad
t\delta_4=(t-k-1)\delta_1 \quad
s\gamma_{3}=({s-k-1})\gamma_{1}.\end{array}
$$
$$\begin{array}{l}
s\gamma_{2}=(k+1)\gamma_{1},\quad
-t\delta_2=(k+1)\delta_1.\end{array}
$$
As before, we prove that a such 1-cocycle defined on
$\mathfrak{sl}(2)$ can be extended to $\mathfrak{osp}(1|2)$. Thus,
we have $\mathrm{dim}\mathrm{H}^1(\frak{osp}
({1|2}),(\frak{D}_{\lambda,\nu;\mu}))=6$.\\

The case: $\mu-\lambda-\nu=k+{3\over2}$  where
$(\lambda,\nu)=(-\frac{s}{2},-\frac{t}{2})$ with $s,\,t\in\{1,\ldots,k\}$
and $s+t\geq k+2$ (super resonance case with $\mu-\lambda-\nu$ semi integer) can be
treated similarly and  we get the same result.\\

%%%%%%%%%%%%%%%%%%%%%%%%%%%%%%%%%%%%%%%%%%%%%%%%%%%%%%%%%%%%%%%%%%%%%%%%%%%%%%%%%%%%%%%%%%%%%%%%%%%%%%%%%%
\section{Singular cases}\label{s5}
%%%%%%%%%%%%%%%%%%%%%%%%%%%%%%%%%%%%%%%%%%%%%%%%%%%%%%%%%%%%%%%%%%%%%%%%%%%%%%%%%%%%%%%%%%%%%%%%%%%%%%%%%%

Finally, we complete the study of the spaces ${\rm
H}^1_\mathrm{diff}(\frak
{osp}(1|2),\mathfrak{D}_{\lambda,\nu;\mu})$ by considering  the
cases $(\lambda,\nu,\mu)$ which are neither super resonant nor
weakly super resonant. We know that the non vanishing spaces
$\HH^1_\mathrm{diff}(\frak {osp}(1|2),\frak{D}_{\l,\nu;\mu})$ only
can appear if $2(\mu-\l-\nu)+1\in\mathbb{N}$, thus, we consider
the following two situations:

\noindent{\bf A.} Let  $\mu-\lambda-\nu=k+1$ where $k\in\mathbb{N}$
and $(\lambda,\nu)=(-{s\over2}, -{t\over2})$ with
$s,\,t\in\{0,\ldots,k+1\}$. In this case the cohomology space is
even and then we have to consider:
 $$
 \begin{array}{l}(\lambda,\nu,\mu),\quad(\lambda+\frac{1}{2},\nu+{1\over2},\mu),
 \quad(\lambda+\frac{1}{2},\nu,\mu+\frac{1}{2})
 \quad\text{and}\quad(\lambda,\nu+\frac{1}{2},\mu+\frac{1}{2}).\end{array}
 $$
 The cases for which  $(\lambda,\nu,\mu)$ is neither super resonant nor weakly super resonant are:\\
 (i) $s+t=k$, in this case only $(\lambda,\nu,\mu)$ is resonant.\\
 (ii) $ s\in\{0,\ldots,k\}$ and
  $t=k+1$, in this case only $(\lambda,\nu+\frac{1}{2},\mu+\frac{1}{2})$ is resonant.\\
  (iii) $s=k+1$ and
  $t\in\{0,\ldots,k\}$, in this case only $(\lambda,\nu+\frac{1}{2},\mu+\frac{1}{2})$ is resonant.\\

\noindent{\bf B.} Let $\mu-\lambda-\nu=k+{3\over2}$ where
$k+1\in\mathbb{N}$ and $(\lambda,\nu)=(-{s\over2},-{t\over2})$
with $s,\,t\in\{0,\ldots,k+1\}$. In this case the cohomology space
is odd and then we have to consider:
$$\begin{array}{l}(\lambda+\frac{1}{2},\nu,\mu),\quad(\lambda,\nu+\frac{1}{2},\mu),
  \quad(\lambda,\nu,\mu+\frac{1}{2})\quad\text{and}\quad(\lambda+\frac{1}{2},\nu+\frac{1}{2},\mu+\frac{1}{2}).\end{array}
 $$
We have to distinguish the following cases:\\
  (i) $(s,t)=(k+1,0)$, in this case only $(\lambda+\frac{1}{2},\nu,\mu)$ and
  $(\lambda,\nu,\mu+\frac{1}{2})$ are resonant.\\
 (ii) $ (s,t)=(0,k+1)$, in this case only $(\lambda,\nu+\frac{1}{2},\mu)$ and
  $(\lambda,\nu,\mu+\frac{1}{2})$ are resonant.\\
  (iii) $ s+t=k+1$ with $st\neq0$, in this case only $(\lambda+\frac{1}{2},\nu+\frac{1}{2},\mu+\frac{1}{2})$
  is non resonant.\\
(iv) $ s=t=k+1$, in this case only $(\lambda,\nu,\mu+\frac{1}{2})$ and
$(\lambda+\frac{1}{2},\nu+\frac{1}{2},\mu+\frac{1}{2})$ are resonant.\\
(v) $ s=k+1$ and $t\in\{1,\dots,k\}$, in this case only $(\lambda,\nu+\frac{1}{2},\mu)$ is non resonant.\\
(vi) $ t=k+1$ and $s\in\{1,\dots,k\}$, in this case only $(\lambda+\frac{1}{2},\nu,\mu)$ is non resonant.\\
 \begin{thm}\label{singular} Let $(\lambda,\nu,\mu)$ be neither super resonant nor
weakly super resonant. \\
(a) If  $\mu-\lambda-\nu=k+1\in\mathbb{N}^*$ and
$(\lambda,\nu)=(-{s\over2}, -{t\over2})$ with
$s,\,t\in\{0,\ldots,k+1\}$ then
\begin{equation*}
{\rm H}^1_\mathrm{diff}(\frak
{osp}(1|2),\mathfrak{D}_{\lambda,\nu;\mu})\simeq \left\{
\begin{array}{llllll}
\mathbb{R} & \hbox{ if }\quad
 s+t=k, \\[2pt]
\mathbb{R}^2 & \hbox{ if }\quad
 s=k+1\,\text{ or } \,t=k+1 \hbox{ with }\,\, k+2\leq s+t\leq2k+1.
\end{array}
\right.
\end{equation*}
(b) If $\mu-\lambda-\nu-\frac{1}{2}=k+1\in\mathbb{N}$ and
$(\lambda,\nu)=(-{s\over2},-{t\over2})$ with
$s,\,t\in\{0,\ldots,k+1\}$ then
\begin{equation*}
{\rm H}^1_\mathrm{diff}(\frak
{osp}(1|2),\mathfrak{D}_{\lambda,\nu;\mu})\simeq
\left\{\begin{array}{llllll}
\mathbb{R}^5 & \hbox{ if }\quad
 s=k+1\hbox{ or } t=k+1\hbox{ with }\,\, k+2\leq s+t\leq2k+1,\\
 \mathbb{R}^3& \hbox{ if }\quad (s,t)=(0,k+1),\,(k+1,0)\hbox{ with }\,\, k\neq-1,\\
 \mathbb{R}^2& \hbox{ if }\quad
 s=t=k+1\hbox{ or }\,\,s+t=k+1\hbox{ but }\,\,st\neq0,\\
 \mathbb{R}& \hbox{ if }\quad
 k=-1.
\end{array}\right.
\end{equation*}
(c) Otherwise ${\rm H}^1_\mathrm{diff}(\frak
{osp}(1|2),\mathfrak{D}_{\lambda,\nu;\mu})=0$.
\end{thm}
\begin{proofname}. Recall that if $2(\mu-\lambda-\nu)+1\notin\mathbb{N}$ then ${\rm H}^1_\mathrm{diff}(\frak
{osp}(1|2),\mathfrak{D}_{\lambda,\nu;\mu})=0$. Thus, assume that
$2(\mu-\lambda-\nu)+1\in\mathbb{N}$.\\

1) \textbf{Even cases}: $\mu-\lambda-\nu=k+1$ where $k\in\mathbb{N}$.\\

i) $(\lambda,\nu)=(-\frac{s}{2},-\frac{k-s}{2})$ with
$s\in\{0,\ldots,k\}$. The restriction of $\Upsilon$ to $\frak
{sl}(2)$ is given by
$$\begin{array}{lll}
\Upsilon(X_{h})(f,g)$=$\alpha_1
h''f^{(s)}g^{(k-s)}+a\displaystyle\sum_{0\leq i\leq
s}\begin{pmatrix}{k+1}\\i\end{pmatrix}h'f^{(i)}g^{(k+1-i)}+
b\displaystyle\sum_{i=s+1}^{k+1}\begin{pmatrix}{k+1}\\i\end{pmatrix}h'f^{(i)}g^{(k+1-i)},
\end{array}$$
$$
\Upsilon(X_{h})(\theta f,\theta
g)=\alpha_2h'f^{(s)}g^{(k-s)},\quad \Upsilon(X_{h})(\theta
f,g)=\theta\alpha_{3}h'f^{(s)}g^{(k-s+1)},$$
$$
\Upsilon(X_{h})( f,\theta g)=\theta\alpha_4h'f^{(s+1)}g^{(k-s)}.$$
As before,  the 1-cocycle condition gives: $$
\Upsilon(X_{h})(F,G)=\alpha_2\left(h'f_{1}^{(s)}g_{1}^{(k-s)}+\theta
h'\left(
f_{0}^{(s+1)}g_{1}^{(k-s)}-f_{1}^{(s)}g_{0}^{(k-s+1)}\right)\right),
$$
that is $a=b=\alpha_1=0$ and $\alpha_{3}=-\alpha_4=-\alpha_2$. We
check that $\Upsilon$ can be extended to $\Pi(\mathcal{H})$ and we
deduce that
 $\mathrm{dim}\HH^1(\frak
{osp}(1|2),(\frak{D}_{\l,\nu;\mu})_1)=1$.\\

ii) $(\lambda,\nu)=(-\frac{s}{2},-\frac{k+1}{2})$ with
$s\in\{1,\ldots,k\}$. The restriction of  $\Upsilon$ to $\frak
{sl}(2)$ is, a priori, given by:
$$\begin{array}{lllll}
\Upsilon(X_{h})(f,g)&=&\alpha_1 \displaystyle\sum_{i=s+1}^{ k+1}
\begin{pmatrix}{k-s}\\i-s-1\end{pmatrix}h'f^{(i)}g^{(k+1-i)},\\
\Upsilon(X_{h})(\theta f,\theta
g)&=&\alpha_2\displaystyle\sum_{i=s}^{
k}\begin{pmatrix}{k-s}\\i-s\end{pmatrix}
h'f^{(i)}g^{(k-i)},\\
\Upsilon(X_{h})(\theta f,g)&=&\alpha_3
\theta\displaystyle\sum_{i=s}^{
k+1}\begin{pmatrix}{k-s+1}\\i-s\end{pmatrix}
h'f^{(i)}g^{(k+1-i)},\\
\Upsilon(X_{h})(f,\theta g)&=&\theta\left(
\alpha_4\displaystyle\sum_{i=s+1}^{
k+1}\frac{k+1}{i}\begin{pmatrix}{k-s}\\i-s-1\end{pmatrix}
h'f^{(i)}g^{(k+1-i)}+\alpha_5h''fg^{(k)} +\alpha_6 h'fg^{(k+1)}
\right).\end{array}$$

The 1-cocycle condition: $ \delta(\Upsilon)(X_x,X_\theta)=0$
gives:
$$
\alpha_6=\alpha_4=0\quad \text{and}\quad
\alpha_1=-\alpha_2=\alpha_3.$$ We easily check that $\Upsilon$ can
be extended to $\Pi(\mathcal{H})$, therefore
 $\mathrm{dim}\HH^1(\frak
{osp}(1|2),\frak{D}_{\l,\nu;\mu})=2$. Of course, we have the same
result if $(\lambda,\,\nu)=(-\frac{k+1}{2},-\frac{s}{2})$ where
$s\in\{1,\ldots,k\}$.\\

2) \textbf{Odd cases}:
$\mu-\lambda-\nu+\frac{1}{2}=k+2\in\mathbb{N}$ and
$(\lambda,\nu)=(-{s\over2},-{t\over2})$ with
$s,\,t\in\{0,\ldots,k+1\}$.\\

i) Let $k=-2$. Here we are in the situation (c) of Theorem
\ref{singular}. Obviously, in this case, $(\lambda,\nu,\mu)$ is
neither super resonant nor weakly super resonant.  The restriction
of $\Upsilon$ to $\frak {sl}(2)$ is, a priori, given by:
\begin{equation*}
\Upsilon(X_h)(F,\,G)=\alpha\theta h'f_0g_0,
\end{equation*}
where $\alpha\in\mathbb{R}$ and $h\in \mathbb{R}_2[x]$. The
1-cocycle condition: $ \delta(\Upsilon)(X_x,X_\theta)=0$ gives the
following equation:
 \begin{equation*}
 x(\Upsilon(X_\theta)(f_0,g_0))'- \Upsilon(X_\theta)(xf'_0,g_0)-
 \Upsilon(X_\theta)(f_0,xg'_0)+{1\over2}\alpha f_0g_0=0.
\end{equation*}
 Thus, we have $\alpha=0$ since ${1\over2}\alpha f_0g_0$ is the unique term in $f_0g_0$
 in the previous equation. By Proposition \ref{sa},
we deduce that, in this case, $\HH^1_\mathrm{diff}(\frak
{osp}(1|2),\frak{D}_{\l,\nu;\mu})=0$.\\

ii) Let  $k=-1$. In this case, $(\lambda,\nu,\mu)$ is neither
super resonant nor weakly super resonant if and only if
$(\lambda,\nu)=(0,0)$. So, let $(\lambda,\nu)=(0,0)$ and consider
the restriction of $\Upsilon$ to $\frak {sl}(2)$  which is, a
priori, given by:
\begin{equation*}
\Upsilon(X_h)(F,G)=h'[\alpha_1f_1g_0+\alpha_2f_0g_1+
\alpha_3\theta f_0g'_0+\alpha_4\theta f'_0g_0+\alpha_5\theta
f_1g_1]+\alpha_6\theta h''f_0g_0.
\end{equation*}
The 1-cocycle condition: $ \delta(\Upsilon)(X_x,X_\theta)=0$, gives
$$\alpha_{1}=\alpha_{2}=-\alpha_{3}=-\alpha_{4}=-\alpha_{6},\quad\alpha_{5}=0.$$
The restriction of  $\Upsilon$ to $\Pi(\mathcal{H})$ can be given by
$$\Upsilon(X_{h_1\theta})(F,G)=\theta\alpha_4h'_1FG.$$
Thus,
$\mathrm{dim}\HH^1_\mathrm{diff}(\frak
{osp}(1|2),\frak{D}_{\l,\nu;\mu})=1$. This proves the situation (b) when $k=-1$.\\

ii) Let $(\lambda,\nu)=(-\frac{k+1}{2},-\frac{k+1}{2})\text{ where
}
 k\in\mathbb{N}^*$.  The restriction of $\Upsilon$ to $\frak {sl}(2)$ is given by
\begin{equation*}\begin{array}{llll}
\Upsilon(X_h)(F,G)&=&\theta(\alpha_1h''f_0g_0^{(k+1)}+
\alpha_2h'f_0g_0^{(k+2)}+\alpha_3h'f_0^{(k+2)}g_0+\beta_1h''f_1g_1^{(k)}\\&~&+
\beta_2h'f_1g_1^{(k+1)}+\beta_3h'f_1^{(k+1)}g_1)+\gamma
h'f_1^{(k+1)}g_0+ \delta h'f_0g_1^{(k+1)}.\end{array}
\end{equation*}
where
$\alpha_1,\,\alpha_2,\,\alpha_3,\,\beta_1,\,\beta_2,\,\beta_3,\,\gamma,\,\delta\in\mathbb{R}$.
From the relation $\delta\Upsilon(X_{h_0},X_{h_1\theta})(F,G)=0$
we deduce that:
\begin{equation*}\begin{array}{llll}
\Upsilon(X_h)(F,G)&=&\theta(\alpha_1h''f_0g_0^{(k+1)}+\beta_1h''f_1g_1^{(k)}).\end{array}
\end{equation*}
We check that $\Upsilon$ can be extended to $\Pi(\mathcal{H})$,
therefore $\mathrm{dim}\HH^1_\mathrm{diff}(\frak
{osp}(1|2),\frak{D}_{\l,\nu;\mu})=2$.

iii) Let $(\lambda,\nu)=(-\frac{k+1}{2},-\frac{t}{2})$ where
$k\in\mathbb{N}^*$ and $t\in\{1,\dots,k\}$. In this case the map
$\Upsilon|_{\mathfrak{sl}(2)}$ can be decomposed as follows:
$\Upsilon|_{\mathfrak{sl}(2)}=B+C+D$ where
\begin{equation*}
B(X_h,F,G)=h''(\theta\beta_1f_0^{(k-t+1)} g_0^{(t)}+\theta \beta_2f_1^{(k-t+1)}
g_1^{(t-1)}+\beta_3f_1^{(k-t)} g_0^{(t)}),
\end{equation*}\begin{equation*}
C(X_h,F,G)=\theta\gamma_1\mathfrak{c}_1(X_h,f_0,g_0)+
\theta\gamma_2\mathfrak{c}_2(X_h,f_1,g_1)+\gamma_3\mathfrak{c}_3(X_h,f_1,g_0),
\end{equation*}\begin{equation*}
D(X_h,F,G)=\delta_1\theta\mathfrak{d}_1(X_h,f_0,g_0)+
\delta_2\theta\mathfrak{d}_2(X_h,f_1,g_1)+\delta_3\mathfrak{d}_3(X_h,f_1,g_0)+
\alpha_1\mathfrak{a}_1(X_h,f_0,g_1)
\end{equation*}where the $\mathfrak{c}_i$, the $\mathfrak{d}_i$
are as those defined in (\ref{c2}) and (\ref{d2}) and $\mathfrak{a}_1$ is as in (\ref{a1}). By the 1-cocycle
relation: $\delta\Upsilon(X_{h_0},X_{h_1\theta})(F,G)=0$ we prove
$$
\delta_1=-\delta_2=-\delta_3,\quad
\gamma_3=\frac{(k+1-t)}{k+1}\,\alpha_{1},\quad
\gamma_2=\frac{-t}{k+1}\,\alpha_{1},\quad
\gamma_1=-\alpha_1.
$$
We prove also that $\Upsilon$ can be extended to
$\Pi(\mathcal{H})$. Thus, $\mathrm{dim}\HH^1(\frak
{osp}(1|2),\frak{D}_{\l,\nu;\mu})_0=5$. Similarly, we study the
case $(\lambda,\nu)=(-\frac{s}{2},-\frac{k+1}{2})$ where
$k\in\mathbb{N}^*$ and $s\in\{1,\dots,k\}$.\\

iv) Let $(\lambda,\nu)=(-\frac{s}{2},-\frac{k+1-s}{2})\text{ where
}
 k\in\mathbb{N}^*\text{ and }s\in\{1,\dots,k\}$.
 In this case the map
$\Upsilon|_{\mathfrak{sl}(2)}$ can be decomposed as follows:
$\Upsilon|_{\mathfrak{sl}(2)}=B+C+D$ where
\begin{equation*}
B(X_h,F,G)=h''(\beta_1f_1^{(s-1)} g_0^{(k-s+1)}+ \beta_2f_0^{(s)}
g_1^{(k-s)}+\theta\beta_4f_0^{(k-t)} g_0^{(t)}),
\end{equation*}\begin{equation*}
C(X_h,F,G)=\theta\alpha_4^0\mathfrak{c}_1(X_h,f_0,g_0)+
\alpha_0^1\mathfrak{c}_3(X_h,f_1,g_0)+\alpha_0^12\mathfrak{c}_4(X_h,f_0,g_1),
\end{equation*}\begin{equation*}
D(X_h,F,G)=\alpha_1^4\theta\mathfrak{d}_1(X_h,f_0,g_0)+\alpha_1^1\mathfrak{d}_3(X_h,f_1,g_0)+
\alpha_1^2\mathfrak{d}_1(X_h,f_0,g_1)+\alpha_3\theta\mathfrak{a}_1(X_h,f_1,g_1)
\end{equation*}where the $\mathfrak{c}_i$, the $\mathfrak{d}_i$
are as those defined in (\ref{c2}) and (\ref{d2}) and $\mathfrak{a}_3$ is as in (\ref{a1}). By the 1-cocycle
relation: $\delta\Upsilon(X_{h_0},X_{h_1\theta})(F,G)=0$ we prove:

$$\begin{cases}
(k+2)\beta_2 =-(k-s+1)\beta_4\\
(k+3)\beta_1=-s\beta_4\\
\end{cases}$$
We prove also that $\Upsilon$ can be extended to
$\Pi(\mathcal{H})$. Thus, $\mathrm{dim}\HH^1(\frak
{osp}(1|2),\frak{D}_{\l,\nu;\mu})_0=2$.\\

\hfill$\Box$\end{proofname}

%%%%%%%%%%%%%%%%%%%%%%%%%%%%%%%%%%%%%%%%%%%%%%%%%%%%%%%%%%%%%%%%%%%%%%%%%%
%%%%%%%%%%%%%%%%%%%%%%%%%%%%%%%%%%%%%%%%%%%%%%%%%%%%%%%%%%%%%%%%%%%%%%%%%%


\begin{thebibliography}{99}
%%%%%%%%%%%%%%%%%%%%%%%%%%%%%%%%%%%%%%%%%%%%%%%%%%%%%%%%%%%%%%%%%%%%%%%%%%%%

\bibitem{ab}
Agrebaoui B and Ben Fraj N, On the cohomology of the Lie
superalgebra of contact vector fields on $S^{1|1}$, {\it Bell.
Soc. Roy. Sci. Li\`ege} {\bf 72}, 6, 2004, 365--375.

\bibitem{b}
S. Bouarroudj, {\it Cohomology of the vector fields Lie algebras
on $\mathbb{R}\mathbb{P}^1$ acting on bilinear differential
operators}, International Journal of Geometric Methods in Modern Physics
 (2005), {\bf 2}; N 1,  23-40.

\bibitem{Fu}
Fuchs  D B, {\it Cohomology of infinite-dimensional Lie algebras},
Plenum Publ. New York, 1986.


\bibitem{gmo} H. Gargoubi, N. Mellouli
and V. Ovsienko {\it Differential Operators on Supercircle:
Conformally Equivariant Quantization and Symbol Calculus,} Letters
in Mathematical Physics (2007) {\bf 79}:51–65.

\end{thebibliography}
\end{document}